\documentstyle[10pt]{amsart}

\newtheorem{prop}{Proposition}
\newtheorem{theo}{Theorem}
\newtheorem{Lemma}{Lemma}
\newtheorem{cor}{Corollary}

\newtheorem{rem}{Remark}

\newcommand{\na}{\nabla}

\newcommand{\Om}{\Omega}

\newcommand{\La}{\Lambda}

\newcommand{\ka}{K{\"a}hler }

\newtheorem{conjecture*}{Conjecture}
\newtheorem{theorem*}{Theorem}
\newtheorem{question*}{Question}

\newcommand{\leftr}{[\hbox{\hspace{-0.15em}}[}
\newcommand{\rightr}{]\hbox{\hspace{-0.15em}}]}

\newcommand{\llra}{\!\!\joinrel{\hbox to 30pt{\rightarrowfill}}}
\newcommand{\lllra}{\!\!\joinrel{\hbox to 50pt{\rightarrowfill}}}
\newcommand{\llllra}{\!\!\joinrel{\hbox to 60pt{\rightarrowfill}}}
\newcommand{\lllllra}{\!\!\joinrel{\hbox to 70pt{\rightarrowfill}}}
\newcommand{\llllllra}{\!\!\joinrel{\hbox to 75pt{\rightarrowfill}}}

\title[Local models and integrability of almost k{\"a}hler 4-manifolds]{
Local models and integrability  of certain almost k{\"a}hler
4-manifolds}
\author{VESTISLAV APOSTOLOV, JOHN ARMSTRONG AND TEDI DR\u{A}GHICI}
\thanks{The first author was supported in part by a PAFARC-UQAM grant,
NSERC grant OGP0023879 and by FCAR grant NC-7264. He is also member of
EDGE, Research Training Network HPRN-CT-2000-00101, supported by the
European Human Potential Programme. The first and the third authors were
supported in part by NSF grant INT-9903302}
\address{Vestislav Apostolov \\ D\'epartement de math\'ematiques\\ UQAM
 \\ succursale Centre-ville  c.p. 8888 \\ Montr\'eal \\ H3C 3P8, Canada}
\email{apostolo@@math.uqam.ca}

\address{John Armstrong \\ 21 Alford House \\
Stanhope Rd \\
London \\
N6 5AL \\ UK}
\email{John.Armstrong@@yolus.com}
\address{Tedi Dr\u{a}ghici \\ Department of Mathematics \\ Florida
International University \\
 Miami FL 33199 \\ USA}
\email{draghici@@fiu.edu}

\begin{document}

\begin{abstract}
We classify, up to a local isometry, all non-K{\"a}hler almost
K{\"a}hler 4-manifolds for which the fundamental 2-form is an eigenform of
the Weyl tensor, and whose Ricci tensor is invariant with respect to
the almost complex structure. Equivalently, such
almost K{\"a}hler 4-manifolds satisfy the third curvature condition of
A. Gray. We use our local classification to show that, in the compact
case, the third curvature condition of Gray is equivalent to the
integrability of the corresponding almost complex structure.

\vspace{0.1cm}
\end{abstract}
\maketitle

\section{Introduction}
Let $(M,\Omega)$ be a symplectic manifold of dimension $2n$.  An
almost complex structure $J$ is called {\it compatible} with the
symplectic form $\Om$, if there exists a Riemannian metric $g$
such that
$$ \Om(\cdot, \cdot) = g(J\cdot, \cdot).$$
In this case, the metric $g$ is also called compatible with $\Om$,
while the triple $(g,J,\Om)$ is referred to as an {\it almost \ka
structure} on $M$. If, additionally, the compatible almost complex
structure $J$ is integrable, then $(g, J, \Om)$ is a K{\"a}hler
structure on $M$. Almost K{\"a}hler manifolds for which the almost
complex structure $J$ is not integrable will be called {\it
strictly} almost K{\"a}hler.

Gromov's theory of pseudo-holomorphic curves \cite{gromov},
Taubes' characterization of Seiberg-Witten invariants of
symplectic 4-manifolds \cite{taubes1,taubes2} and the recent works
of Donaldson \cite{donaldson1,donaldson2} conclusively situated
the study of (strictly) almost \ka manifolds into the mathematical
main-stream. A very recent work of LeBrun \cite{lebrun3}, for
instance, inspired by some aspects of Taubes' construction of
solutions of the Seiberg-Witten equations on symplectic
4-manifolds, relates the existence of strictly almost \ka
4-manifolds having particular properties of the curvature with
some fundamental problems in Riemannian Geometry, such as the
existence and uniqueness of Einstein metrics, or of metrics
minimizing volume among all Riemannian metrics satisfying
point-wise lower bounds on sectional curvatures.

Despite the now compelling appeal of the ``Riemannian'' aspect of
almost \ka geometry, the subject suffered for a long time from a
genuine lack of interesting examples. For instance, an old still
open conjecture of Goldberg \cite{Go} affirms that there are no
{\it Einstein}, strictly almost \ka metrics on a compact
symplectic manifold. K. Sekigawa proved \cite{Se2} that the
conjecture is true if the scalar curvature is non-negative. The
case of negative scalar curvature is still wide open, despite of
the recent discovery  of {\it complete} Einstein strictly almost
\ka manifolds of any dimension $2n \ge 6$ \cite{ADM}, and of
(local) Ricci-flat strictly almost \ka metrics of dimension four
\cite{ACG,Arm1,NuP}. Indeed, in addition to some subtle
topological obstructions in dimension four
\cite{kotschick,lebrun1,lebrun2,lebrun3}, it turns out that there
is a number of (rather not well understood yet) local obstructions
to the existence of Einstein strictly almost \ka metrics
compatible with a given symplectic form \cite{AA}. Further
progress on the Goldberg conjecture seems to hinge on a better
understanding of these obstructions.

\vspace{0.2cm} Classically, on a compact manifold $M$, the
Einstein condition appears as the Euler-Lagrange equation of the
{\it Hilbert functional} ${\bf S}$, the total scalar curvature,
acting on the space of all Riemannian metrics on $M$ of a given
volume. A ``symplectic'' setting of this variational problem was
first considered by Blair and Ianus \cite{BI}: let $(M,\Om)$ be a
compact symplectic manifold and restrict the functional ${\bf S}$
to the space of $\Om$-compatible Riemannian metrics, then the
critical points are the almost \ka metrics $(g,J,\Om)$ whose
Ricci tensor ${\rm Ric}$ is $J$-invariant, i.e. satisfies:
\begin{equation}\label{eq1}
{\rm Ric}(J\cdot, J\cdot) = {\rm Ric}(\cdot,\cdot).
\end{equation}
The Euler-Lagrange equation (\ref{eq1}) is a weakening
of both the Einstein and the \ka conditions.  Furthermore, Blair
\cite{Bl} observed that for any almost \ka metric $(g,J)$ the
following relation holds:
$$\frac{1}{4} \int_M |\nabla J|^2 d\mu  + {\bf S}(g) =  \frac{4\pi}{(n-1)!}
(c_1 \cdot [\Omega]^{\wedge (n-1)})(M),$$ where $\nabla$ is the
Riemannian connection of $g$, $|\cdot |$ is the point-wise norm
induced by $g$, and $c_1$ and $d\mu = \frac{1}{n!}\Omega^{\wedge
n}$ are respectively the first Chern class and the the volume
form of $(M,\Om)$. Note that the right hand-side of the above
equality is a symplectic invariant (hence is independent of the
choice of a particular $\Om$-compatible metric); thus the
$\Om$-compatible almost \ka structures satisfying (\ref{eq1}) are
also the critical points of the {\it Energy functional}, which
acts on the space of $\Om$-compatible almost complex structures by
$${\bf E} (J) = \int_M |\na J|^2 d\mu.$$
From this point of view, critical almost \ka metrics have been
recently studied in \cite{Le-Wang}. Clearly, the functional ${\bf
E}$ (resp. ${\bf S}$) is bounded from below (resp. from above),
the \ka metrics being minima of ${\bf E}$ (resp. maxima of ${\bf
S}$). However, a direct variational approach of finding extremal
metrics for these functionals seems not to be easily applicable
as it may happen that the infimum of ${\bf E}$ be zero, although
$M$ does not carry \ka structures at all (see \cite{Le-Wang}). In
fact, apart from the explicit compact examples of \cite{AGI,DM}
(which, multiplied by \ka manifolds, provide examples of any
dimension $2n\ge 6$), almost nothing seems to be known in general
about existence of compact strictly almost \ka metrics satisfying
(\ref{eq1}). In particular, no such example is yet known on a
compact symplectic 4-manifold.

One reason for which some technical difficulties appear in
applying global methods is perhaps hidden in the fact that even
locally the equation (\ref{eq1}) is difficult to be solved. Indeed,
it consists of a system of PDE's for a compatible metric with $\Om$, which
does not satisfy the Cartan test \cite{AA}. One finds instead a number of
nontrivial relations between higher jets of $J$ and the curvature of
$g$.  On a compact manifold precisely these relations  are at the origin
of the integrability results obtained in
\cite{AA,ADK,Arm1,Dr2,Se2}. The few known (local) examples
appear, in fact, as a byproduct of additional geometric structure which
makes the study of the Euler-Lagrange equation (\ref{eq1}) more tractable.

\vspace{0.2cm} The additional structure relevant to this paper
comes as an extra-assumption imposed on the curvature tensor of
the almost \ka 4-manifold, namely
\begin{equation}\label{eq2}
W(\Om) = \nu \Omega,
\end{equation}
where  $W$ is the Weyl tensor of $g$, viewed as a symmetric
traceless endomorphism acting on the space of 2-forms; thus,
equation (\ref{eq2}) is equivalent to $\Om$ being an eigenform of
the Weyl tensor, in which case the eigenfunction $\nu$ is smooth
and equal, up to a factor $1/6$, to the so-called conformal
scalar curvature $\kappa$ of $(g,J)$ (see Section 2.3 for a
precise definition).

The condition (\ref{eq2}) in its own right appears to be natural
in the context of four dimensional almost \ka geometry. Indeed,
motivated by harmonic maps theory, C. Wood \cite{wood} showed
that  an almost \ka 4-manifold $(M,g,J,\Om)$ satisfies
(\ref{eq2}) if and only if $J$ is a critical point of the energy
functional ${\bf E}$, but under variations through all almost
complex structures compatible with the given metric $g$; for this
reason such almost complex structures were called {\it harmonic}
in \cite{wood}. Furthermore, almost \ka metrics saturating the
new curvature estimates of \cite{lebrun3} must all satisfy
(\ref{eq2}).

Let us finally mention that each of the conditions (\ref{eq1}) and
(\ref{eq2})  corresponds to the vanishing of a certain
irreducible component of the curvature $R$ under the action of
$U(2)$ on the space of algebraic curvature tensors \cite{TV}.
Taken together, (\ref{eq1}) \& (\ref{eq2}) are equivalent to the
more familiar {\it third curvature condition of Gray} (see
\cite{Gr}):
\begin{equation}\label{gray3}
 R_{XYZU} = R_{JXJYJZJU}.
\end{equation}
Evidently, the curvature of a \ka manifold satisfies
(\ref{gray3}), hence any of the conditions (\ref{eq1}) and
(\ref{eq2}).

\vspace{0.2cm} In \cite{AAD}, we have constructed an explicit
family of strictly almost \ka metrics that all satisfy the third
Gray condition $(\ref{gray3})$.  The main purpose of this paper
is to close the circle of ideas from our previous work
\cite{AAD}, by showing that, conversely, any strictly almost \ka
4-manifold which satisfies $(\ref{gray3})$ is locally modeled by
a metric in this family:

\begin{theo}\label{thm1} Let $(\Sigma, g_{\Sigma}, \Om_{\Sigma})$ be an
oriented Riemann surface with metric $g_{\Sigma}$ and volume form
$\Omega_{\Sigma}$, and let $h = w + iv$ be a non-constant holomorphic
function on $\Sigma$, whose real part $w$ is everywhere positive.
On the product of $\Sigma$ with ${\mathbb R}^2 = \{ (z,t) \}$
consider the symplectic form
\begin{equation}\label{om0}
\Om = \Om_{\Sigma} - dz\wedge dt
\end{equation}
and the compatible Riemannian metric
\begin{equation}\label{g0}
g = g_{\Sigma} + w dz^{\otimes 2} + \frac{1}{w} (dt +
vdz)^{\otimes 2};
\end{equation}
Then,
\begin{enumerate}
\item[\rm (i)] $(g,\Om)$ defines a strictly almost \ka structure whose
curvature satisfies the third Gray condition {\rm (\ref{gray3})}.
\item[\rm (ii)] For any connected strictly almost \ka  4-manifold
$(M,g,\Om)$ whose curvature satisfies {\rm (\ref{gray3})} there
exists an open dense subset $U$ with the following property:  in a
neighborhood of any point of $U$, $(g,\Om)$ is homothetic to an
almost \ka structure given by {\rm (\ref{om0}-\ref{g0})}.
\end{enumerate}
\end{theo}

This local classification of strictly almost \ka 4-manifolds
satisfying (\ref{gray3}) includes as particular cases some
previously known results from \cite{AAD} and \cite{Arm1}. Indeed,
the {\it Einstein} metrics in our family correspond to the case
when $(\Sigma, g_{\Sigma})$ is a flat surface so that $g$ is a
Gibbons-Hawking metric \cite{GH}  with respect to a
translation-invariant harmonic function; we get precisely the
selfdual Ricci-flat strictly almost \ka examples of
\cite{Arm1,NuP,Tod}. Furthermore, if we consider the
half-plane realization of the hyperbolic  space $({\mathbb H}^2,
g = \frac{dx^{2} + dy^{2}}{x^2})$ and take
the holomorphic function
$h= x+ iy$, then (\ref{om0}-\ref{g0}) gives the four dimensional proper
3-symmetric space
$({\rm Isom}({\mathbb E}^2)\cdot Sol_2)/SO(2)$ (see
\cite{gray,kow}), which can be characterized as the only strictly almost
\ka 4-manifold whose curvature tensor preserves the type decomposition of
complex 2-forms
\cite{AAD} (a condition known as the {\it second curvature condition of
Gray}
\cite{Gr}).

To prove Theorem \ref{thm1} we use local methods which  have
originated from \cite{Arm1}, and have been further developed in
\cite{AAD}. However, the more general result of the current paper
requires two main new ingredients.  The first comes in the form
of a Unique Continuation Principle for the Nijenhuis tensor, $N$,
of an almost \ka 4-manifold satisfying (\ref{gray3}) (Proposition
\ref{ucp}), which allows us to realize the dense set $U$ in
Theorem \ref{thm1} as the  subset of points where $N$ does not
vanish. We then consider on $U$ the orthogonal splitting of the
tangent bundle $TM$ into the sum of two linear complex
sub-bundles, $ D$ and ${D}^{\perp}$, where ${D} = \{ TM \ni X :
\na_X J =0 \}$ is the so-called {\it \ka nullity} of $J$, while
its orthogonal complement, ${D}^{\perp}$, is identified with the
canonical bundle $K_J$ of $J$ (see Sec. 2.2). This allows us to
define another almost complex structure, $I$, which is equal to
$J$ when acting on ${ D}$, but to $-J$ on ${D}^{\perp}$. The
almost complex structure $I$ is compatible with $g$, and induces
the opposite orientation on $M$. Motivated by this fact, we shall
often refer to $(g,I)$ as the {\it opposite} almost Hermitian
structure associated to $(g,J)$. The second new ingredient, which
is the heart of the proof of the above Theorem \ref{thm1},
consists of showing that for any strictly almost \ka metric
$(g,J)$ satisfying (\ref{gray3}), the opposite almost Hermitian structure
$(g,I)$ is, in fact, a \ka structure (Proposition \ref{localadd-prop}).
Once this is achieved, Theorem \ref{thm1} follows by
\cite[Theorem 2]{AAD} (see also Theorem \ref{thm3} below).

\vspace{0.2cm} Note that the metric (\ref{g0}) of Theorem
\ref{thm1} is endowed with an ${\mathbb R}^2$-isometric action
which is {\it surface-orthogonal} and preserves the symplectic
form (\ref{om0}); such an action cannot be extended to a
Hamiltonian toric action on a {\it compact} symplectic 4-manifold.
One is then lead to suspect that there will be no compact examples
of strictly almost \ka 4-manifolds satisfying (\ref{gray3})
(although there are complete ones \cite{AAD}). We confirm this
with the following integrability theorem, which is a direct
generalization of results of \cite{ADK} and \cite{Arm1}:

\begin{theo}\label{thm2} A compact almost
\ka 4-manifold whose curvature satisfies the third Gray condition
{\rm (\ref{gray3})} is necessarily \ka.
\end{theo}

The proof of Theorem \ref{thm2} involves the local result of
Theorem \ref{thm1} and the unique continuation property
established in Proposition \ref{ucp}, together with some global
arguments relying on the Enriques-Kodaira classification of
compact complex surfaces and structure results for non-singular
holomorphic foliations \cite{brunella,mont} and authomorphism
groups \cite{maruyama} of ruled surfaces.

\vspace{0.2cm}

This article can be seen as a natural continuation (and
conclusion) of our previous work \cite{AAD}. The authors have
nonetheless endeavored to make the current paper as
self-contained as possible. The necessary background of almost
\ka geometry and a quick review of previous results is displayed
in Sections 2 and 3 below. The proofs of Theorems \ref{thm1} and
\ref{thm2} are then presented in Sections 4 and 5.

\section{Elements of almost \ka geometry}

\subsection{Type-decompositions of forms and vectors}
Throughout the paper, $(M, g, J, \Om)$ will denote an almost
K{\"a}hler manifold  of (real) dimension $4$, where: $g$ is a
Riemannian metric, $J$ is a $g$-orthogonal almost-complex
structure  -- i.e. $(M,g,J)$ is an almost Hermitian manifold -- and
$\Om(\cdot,
\cdot)= g(J\cdot,
\cdot)$ is the fundamental 2-form of $(g,J)$,  which is closed, and
therefore $(M,\Om)$ is a symplectic manifold.

We denote by:  $TM$ the (real) tangent bundle of $M$; $T^*M$ the
(real) cotangent bundle; $\La^rM, r=1,...,4$ the bundle of real
$r$-forms; $\langle \cdot,\cdot \rangle$ the inner product
induced by $g$ on these bundles (or on their tensor products),
with the following conventions for the wedge product and for the
inner product on forms:
$$(\alpha_1\wedge...\wedge \alpha_r)_{X_1,...,X_r} = {\rm
det}((\alpha_i(X_j)),$$
$$\langle \alpha_1\wedge...\wedge \alpha_r, \beta_1\wedge...\wedge \beta_r
\rangle = {\rm det}(\langle \alpha_i, \beta_j \rangle),$$
where ${\rm det}$ denotes the
determinant.

Using the metric,  we shall implicitly identify vectors and
covectors and, accordingly,  a 2-form $\phi$ with the
corresponding skew-symmetric endomorphism of the tangent bundle
$TM$, by putting: $\langle \phi(X),Y \rangle = \phi(X,Y)$ for any vector
fields $X,Y$. Also, if $\phi, \psi \in T^*M^{\otimes 2}$, by $
\phi \circ \psi$ we understand the endomorphism of $TM$ obtained
by the composition of the endomorphisms corresponding to the two
tensors.

The almost-complex structure $J$ gives rise to a type
decomposition of complex vectors and forms. By convention, $J$
acts on the cotangent bundle $T^*M$ by $(J\alpha)_X =
-\alpha_{JX}$, so that $J$ commutes with the Riemannian duality
between $TM$ and $T^*M$. We shall use the standard decomposition
of the complexified cotangent bundle
$$T^*M\otimes {\Bbb C}= \La^{1,0}M \oplus \La^{0,1}M,$$
given by the $(\pm i)$-eigenspaces of $J$, the type decomposition of
complex
$2$-forms
\begin{equation*}\label{Lambda2c}
\La^2M\otimes {\mathbb C} = \La^{1,1}M \oplus \La^{2,0}M \oplus
\La^{0,2}M,
\end{equation*}
and the type decomposition of symmetric (complex) bilinear forms
$$S^2M \otimes {\mathbb C} = S^{1,1}M \oplus S^{2,0}M \oplus S^{0,2}M.$$

Besides the type decomposition of complex forms, we shall also
consider the splitting of real 2-forms into $J$-invariant and
$J$-anti-invariant ones; for any  real section of $\La^2M$ (resp.
of $S^2 M$), we shall use the super-script $'$ to denote the
projection to the real sub-bundle $\La^{1,1}_{\Bbb R} M$ (resp.
$S^{1,1}_{\Bbb R} M$) of $J$-invariant 2-forms (resp. symmetric
2-tensors), while the super-script $''$ stands for the projection
to the bundle $\leftr \La^{0,2}M \rightr $ (resp. $\leftr S^{0,2}M
\rightr$) of $J$-anti-invariant ones; here and henceforth $\leftr
\  \cdot \ \rightr$ denotes the real vector bundle underlying  a
given complex bundle. Thus, for any section $\psi$ of $\La^2 M$
(resp. of $S^2 M$), we have the orthogonal splitting $\psi = \psi'
+ \psi'' ,$ where
$$\psi'(\cdot , \cdot) = \frac{1}{2}(\psi(\cdot , \cdot) + \psi(J\cdot,
J\cdot)) \; \mbox{ and } \;
\psi''(\cdot , \cdot) = \frac{1}{2}(\psi(\cdot , \cdot) - \psi(J\cdot,
J\cdot)) .$$ We finally define the $U(2)$-decomposition (with respect to
$J$) of real 2-forms
\begin{equation}\label{Lambda2r}
\La^2M = {\mathbb R}\cdot \Om \oplus \La^{1,1}_0 M \oplus \leftr
\La^{0,2}M
\rightr,
\end{equation}
where $\La^{1,1}_0 M$ is the sub-bundle of the {\it primitive}
(1,1)-forms, i.e. the $J$-invariant 2-forms which are point-wise
orthogonal to $\Om$. The above splitting fits in with the well
known $SO(4)$-decomposition of $\La^2M$
\begin{equation*}\label{SO(4)La2}
\La^2 M = \La^+ M \oplus \La^-M
\end{equation*}
into the sub-bundles $\La^{\pm}M$ of selfdual (resp.
anti-selfdual) forms. Indeed, we have
\begin{equation}\label{U(2)La+}
\La^+M = {\mathbb R}\cdot \Om \oplus \leftr
\La^{0,2}M \rightr, \ \ \
\La^-M = \La^{1,1}_0 M.
\end{equation}
Note that $S^{1,1}_{\Bbb R}M$ can  be identified with
$\La^{1,1}_{\Bbb R}M$ via the complex structure $J$: for any $S
\in S^{1,1}_{\Bbb R} M$,
$$(S\circ J)(\cdot,\cdot) : = S(J\cdot, \cdot) $$
is the corresponding element of $\La^{1,1}_{\Bbb R}M$. Also, the
real bundle $\leftr \La^{0,2}M \rightr$ (resp.  $\leftr
S^{0,2}M\rightr $) inherits a canonical complex structure, still
denoted by ${J}$, which is given  by
$$({J}\psi)(X,Y) := -\psi(JX,Y),
\ \forall \psi \in \leftr \La^{0,2}M \rightr, $$ so that $(\leftr
\La^{0,2}M \rightr, J)$ becomes isomorphic to the {\it
anti-canonical} bundle $K^{-1}_J(M)$ of $M$. We adopt a similar
definition for the action of ${J}$ on $\leftr S^{0,2}M\rightr $.
Notice that, using the metric $g$, $\leftr S^{0,2}M\rightr$ can
be also viewed as the bundle of symmetric, $J$-anti-commuting
endomorphisms of $TM$.

\subsection{The ${U}(2)$-decomposition of the curvature}

The Riemannian curvature $R$ is defined by $R _{X, Y} Z = (\na
_{[X, Y]} - [\na_X, \na_Y])Z,$ where $\na$ is the Levi-Civita
connection of $g$. Using the metric, it will be considered as a
section of the bundle $S^2 (\Lambda ^2M)$ of symmetric
endomorphisms of $\Lambda ^2 M$, or of the tensor product $\La^2
M \otimes \La^2M$, depending on the context. The conformal part
of $R$, the Weyl tensor $W$, commutes with the Hodge operator $*$
acting on 2-forms and, accordingly,  splits as $W = W^+ + W^-$,
where $W^{\pm} = \frac{1}{2} (W \pm W \circ *)$. $W^+$ is called
the {\it selfdual Weyl tensor}; it acts trivially on $\Lambda ^-
M$ and will be considered in the sequel as a field of (symmetric,
trace-free) endomorphisms of $\Lambda ^+ M$. Similarly, the {\it
anti-selfdual Weyl tensor}, $W^-$, will be considered as a field
of endomorphisms of $\Lambda ^ - M$.

The Ricci tensor, ${\rm Ric}$, is the symmetric bilinear form
defined by ${\rm Ric} (X, Y) = {\rm tr} \, \{ Z \to R _{X, Z}
Y\}$; alternatively, ${\rm Ric} (X, Y) = \sum _{i = 1}^4 \langle
R _{X, e _i} Y, e _i \rangle$ for any $g$-orthonormal basis $\{ e
_i \}$. We then have ${\rm Ric} = \frac{s}{4} \, g + {\rm Ric}
_0$, where $s$ is the scalar curvature (that is, the trace of
${\rm Ric}$ with respect to $g$)  and ${\rm Ric} _0$ is the {\it
trace-free Ricci tensor}. The latter can be made into a section
of $S ^2 (\Lambda ^2 M)$, then denoted by $\widetilde{{\rm Ric}
_0}$, by putting $ \widetilde{{\rm Ric} _0} (X \wedge Y) = {\rm
Ric} _0 (X) \wedge Y + X \wedge {\rm Ric} _0 (Y);$ equivalently,
for any section $\phi$ of $\La^2M$ we have
$$\widetilde{{\rm Ric} _0} (\phi) = {\rm Ric}_0\circ \phi +
\phi \circ {\rm Ric}_0.$$
It is   readily checked  that  $\widetilde{{\rm Ric} _0}$ satisfies the
first Bianchi
identity, i.e. $\widetilde{{\rm Ric} _0}$ is a tensor of the same kind as
$R$ itself,  as well as $W
^+$ and $W^-$; moreover, $ \widetilde{{\rm Ric} _0}$
anti-commutes with
$*$, so that it can be viewed as a field of homomorphisms from $\Lambda ^+
M$ into $\Lambda ^- M$, or from $\Lambda ^-
M$ into $\Lambda ^+ M$ (adjoint to each other);
we eventually get  the well-known Singer-Thorpe decomposition  of $R$, see
e.g. \cite{besse}:
\begin{equation*}\label{SO(4)R} R = \left(
\begin{array}{c|c} W^+ +
\frac{s}{12}
\, {\rm
  Id} _{| \Lambda ^+ M}   & (\frac{1}{2}\widetilde{{\rm Ric} _0})
_{|
\Lambda ^- M} \\
  \\ \hline \\ (\frac{1}{2}\widetilde{{\rm Ric} _0})_{| \Lambda ^+ M}
 & W ^- +
\frac{s}{12} \, {\rm
  Id} _{| \Lambda ^- M}
\end{array} \right).
\end{equation*}

\vspace{0.2cm}

In the presence of a $g$-orthogonal almost complex structure $J$,
which induces the chosen orientation of $M$,  the above
$SO(4)$-decomposition can be further refined to get seven
irreducible $U(2)$-invariant pieces \cite{TV}. To see this, we
first decompose the traceless Ricci tensor into its $J$-invariant
and $J$-anti-invariant parts, ${\rm Ric}_0'$ and ${\rm
Ric}_0''={\rm Ric}''$, which gives rise to  the decomposition
\begin{equation*}\label{U(2)ric}
\widetilde{{\rm Ric}_0} = \left(
\begin{array}{c|c} \widetilde{{\rm Ric}_0^{'}} &
\widetilde{{\rm Ric}_0^{''}} \end{array} \right) .
\end{equation*}
We shall denote by $\rho$ the {\it Ricci form} of $(M,g,J)$,
which is the $(1,1)$-form corresponding to ${\rm Ric}'$ via $J$
(we thus have  $\rho = {\rm Ric}' \circ J$), and by $\rho_0 =
{\rm Ric}_0' \circ J$ its primitive part.

Furthermore, with respect to (\ref{U(2)La+}), $W^+$ decomposes as
\begin{equation*}\label{U(2)W}
W^+= \left( \begin{array}{c|c}   \frac{\kappa}{6}
      &  W^+_2 \\ \\
\hline \\
                          (W^+_2)^*   &
-\frac{\kappa}{12}{\rm Id}|_{\leftr \La^{0,2} M \rightr} + W^+_3
\end{array}
\right),
\end{equation*}
where:
\begin{enumerate}
\item[$\bullet$] the smooth function $\kappa $ is the so called {\it
conformal scalar curvature} of $(g, J)$, defined by $\kappa =
3\langle W^+(\Om), \Om \rangle$; thus, $\kappa$ is determined by the scale
action of $W^+$ on any of the  two factor of (\ref{U(2)La+}).

\item[$\bullet$] the component $W^+_2$ is the piece of $W^+$ that
interchanges the two factors of (\ref{U(2)La+}); it  thus
can be seen  as a section of
$K^{-1}_J(M)$ by the following equality
\begin{equation}\label{w^+2}
W^+_2=\frac{1}{2}\Big(\rho_*'' \otimes \Omega + \Omega\otimes
\rho_*''\Big),
\end{equation}
where $\rho_*''$ is the $J$-anti-invariant component of the so-called
{\it star-Ricci form} $\rho_* = R(\Om)$ of $(g,J)$; recall that the
{\it star-Ricci tensor} ${\rm Ric^*}$ is then defined by
${\rm Ric^*}(\cdot,
\cdot) = -
\rho_*(J\cdot, \cdot)$;  in particular,
$\rho_*''=0$ if and only if $W^+(\Om) = \frac{\kappa}{6}\Om$;

\item[$\bullet$] the component $ W^+_3$ can be viewed as a section
of $K_J^{-2}(M)$; with respect to  any (local)
section
$\phi $ of
$\leftr
\La^{0,2}M
\rightr$, it can be written as
\begin{equation}\label{w^+3}
W^+_3 = \frac{\lambda}{2}[\phi\otimes \phi - J\phi\otimes J\phi]
      + \frac{\mu}{2}[\phi\otimes J\phi + J\phi\otimes\phi],
\end{equation}
$\lambda $ and $\mu$ being (locally defined) smooth functions;
equivalently,
$W^+_3$ is the component of $W^+$ that belongs to the bundle of symmetric
endomorphisms of
$\leftr
\La^{0,2}M
\rightr$ which anti-commute with $J$.
\end{enumerate}
We thus get the $U(2)$-splitting of the Riemannian curvature
of an almost Hermitian 4-manifold \cite{TV}:
\begin{equation}\label{u(2)}
R = \frac{s}{12} {\rm Id}_{| \La^2M} + W_1 ^+ + W_2 ^+ + W_3 ^+
 + {\widetilde{{\rm Ric}'_0}}
 + {\widetilde {{\rm Ric}''_0}} + W^- ,
\end{equation}
where $W^+_1$ denotes the scalar component of $W^+$, specifically
given by
\begin{equation}\label{w^+1}
W_1^+ = \frac{\kappa}{8} \Om \otimes \Om - \frac{\kappa}{12} {\rm
Id}_{| \La^+M}.
\end{equation}

The following immediate application  of the  decomposition
(\ref{u(2)}) provides alternative definitions (in dimension 4)
for the third curvature condition of Gray (\ref{gray3}):
\begin{Lemma}\cite{TV}\label{lem1}
For an almost Hermitian 4-manifold $(M, g, J)$ the following
conditions are equivalent:
\begin{enumerate}
\item[\rm (i)] ${\rm Ric}''=0$ and $\rho_*'' =0$.
\item[\rm(ii)]
${\rm Ric}''=0$ and  $W(\Om) = \frac{\kappa}{6}\Om$.
\item[\rm (iii)] $ {\rm Ric^*} - {\rm Ric} = \frac{(\kappa - s)}{6} g$.
\item[\rm (iv)] $ R_{XYZU} = R_{JXJYJZJU}$.
\item[\rm (v)] $R$ preserves the type decomposition {\rm (}with respect
to
$J${\rm )} of real 2-forms.
\end{enumerate}
\end{Lemma}

\subsection{The Nijenhuis tensor and the curvature}
A central object in our study is the {\it Nijenhuis tensor} (or
{\it complex torsion}) of an almost complex structure $J$,
defined by
$$N_{X,Y} = [JX,JY] - J[JX,Y] - J[X,JY] - [X,Y].$$
Thus, $N$ is a 2-form with values in $TM$, which by
Newlander-Nirenberg theorem \cite{Ne-Ni} vanishes if and only if
$J$ is integrable.

Alternatively, the Nijenhuis tensor can be viewed as a  map from
$\La^{1,0}M$ to  $\La^{0,2}M$ in the following manner: given a
complex (1,0)-form $\psi$, we define by $\partial \psi$ and
${\bar \partial} \psi$ the projectors of $d\psi$ to  $\La^{2,0}M$
and $\La^{1,1}M$, respectively. In general, $d\neq
\partial + {\bar \partial}$, as $d\psi$ can also have a component of type
(0,2), which we denote by $N(\psi)$; writing
$\psi = \alpha + iJ\alpha$ where $\alpha$ is a real 1-form,  one
calculates
$$N(\psi)_{X,Y}= \frac{1}{2}\alpha(N_{X,Y})=\frac{1}{4}
\psi(N_{X,Y}), \ \ \forall X,Y \in T^{0,1}M.$$
Equivalently, for every real 1-form $\alpha$ we have
\begin{equation}\label{dc}
(d^J\alpha)_{X,Y}'' = (J(d\alpha)'')_{X,Y} +
\frac{1}{2}J\alpha (N_{X,Y}),
\end{equation}
where  $d^J = i({\bar {\partial}} - \partial)$; thus, $d^J$ acts
on real functions by $d^J f = Jdf,$  on  real 1-forms by $(d^J
\alpha)_{X,Y} = -d(J\alpha)_{JX,JY},$ etc. Observe  that $d^2=0$
does not imply $d d^J  + d^J d=0$, unless $N\equiv 0$. For
example, using (\ref{dc}) one obtains on functions
\begin{equation} \label{ddJ+dJd}
 (d d^J + d^J d) f = 2(dd^J f)'' =
(d^Jf)_{N(\cdot, \cdot)}.
\end{equation}
In particular, $dd^J f$ is not, in general, a $(1,1)$-form;
however, for an almost \ka manifold its projection to the
$\Om$-factor of (\ref{Lambda2r}) is a  second order linear
differential operator which, up to a sign, is nothing else than
the Riemannian Laplacian of $(M,g)$:
\begin{equation}\label{laplacian}
\Delta f =  -\langle dd^J f ,\Om \rangle.
\end{equation}

 For an almost \ka manifold $(M,g,J,\Om)$, the
vanishing of $N$ is equivalent to $\Om$ being parallel
with respect to the Levi-Civita connection $\na$ of $g$;
specifically, for any almost \ka manifolds the following
identity  holds (cf. {\it e.g.} \cite{KN2}):
\begin{equation} \label{naOm-N}
(\na_X\Om)(\cdot,\cdot) = \frac{1}{2} \langle JX, N(\cdot,\cdot)
\rangle.
\end{equation}
Using the above relation, we shall often think of $N$ as a
$T^*M$-valued 2-form, or as a $\Lambda^2 M$-valued 1-form,
by tacitly identifying $N$ with $\na \Om$.
Further, since $\Om$ is closed and $N$ is a $J$-anti-invariant
2-form with values in $TM$, one easily deduces that $\na \Om$ is,
in fact, a section of the vector bundle $\leftr \La^{0,1}M\otimes
\La^{0,2}M\rightr$, i.e. the following relation holds:
\begin{equation}\label{naJ}
\na_{JX} J = -J(\na_X J), \ \ \forall X \in TM.
\end{equation}
Note that, in four dimensions, any section of $\leftr \La^{0,1}M\otimes
\La^{0,2}M
\rightr$ has a non-trivial kernel; then,  at any point $x\in M$,
we define the sub-space
$$D_x= \{ X\in T_xM : (\na_X J)_x=0 \},$$
which we call the {\it \ka nullity} of $(g,J)$; thus, ${\rm dim
}_{\mathbb R}(D)= 4$ or $2$, depending on whether or not $N$
vanishes at the given point. Assuming that $N \neq 0$, we let
$D^{\perp}$ be the orthogonal complement of $D$ in $TM$. Relation
(\ref{naJ}) identifies the complex line bundle $D^{\perp}$ with
the {\it canonical} bundle $K_J(M) \cong \La^{2,0}M$;
consequently, we have the following relations for the
corresponding first Chern classes:
$$c_1(D^{\perp})= -c_1(TM) = -c_1(J) \;, \; \;
 c_1(D) = c_1(TM) - c_1(D^{\perp}) = 2c_1(J).$$
On the open subset of points where $N \neq 0$, we shall also
consider the opposite almost complex structure $I$, which is equal
to $J$ on $D$, but to $-J$ on $D^{\perp}$; clearly $I$ is
compatible with the metric but yields the opposite orientation to
the one of $(M,J)$. We mention below the following  observation
due to J. Armstrong \cite{Arm0} and C. LeBrun \cite{lebrun3}:
\begin{Lemma}\label{lem2}  Let
$(M,g,J,\Om)$ be an almost-\ka 4-manifold for which
the Nijenhuis tensor nowhere vanishes. Then, the Chern classes of
the almost complex structures $J$ and $I$ satisfy $c_1(I)=
3c_1({J})$; if, moreover,  $M$ is compact, then the Euler
characteristic $\chi(M)$ and signature $\sigma(M)$ of $M$ are
related by $5\chi(M) = - 6\sigma(M).$
\end{Lemma}
\noindent{\it Proof.} We clearly have $$c_1(I) = c_1(D)- c_1(D^{\perp})=
3c_1(J).$$ By Wu's formula  we know that $c_1^2(M,J) = 2\chi(M) +
3\sigma(M)$ and $c_1^2(M,I)= 2\chi(M) - 3\sigma(M)$, so that we
get $-9(2\chi(M) + 3\sigma(M))= 2\chi(M) - 3\sigma(M)$; the claim
follows. $\Box$

\vspace{0.2cm}

For a \ka structure $(g,J,\Om)$ some of the curvature components
defined in Section 2.2 identically vanish: indeed, in this case
${\rm Ric}''=0$, $W^+_2=0$, $W^+_3=0$ and the scalar function
$\kappa$ is nothing else than the scalar curvature $s$. This
suggests that for an almost \ka structure those components depend
directly on the Nijenhuis tensor (and its covariant derivative).
To see this, consider the Ricci identity for $\Om$:
\begin{equation}\label{ricid}
(\na^2_{X,Y} - \na ^2_{Y,X})(\Om) = [J, R_{X,Y}],
\end{equation}
where $R_{X,Y}$ is viewed as an endomorphism of $TM$ and $[\cdot, \cdot]$
denotes the commutator. A contraction of this relation leads to
\begin{equation} \label{ric-ric*-gen}
-\sum_{i=1}^{4} (\na^2_{e_i,Z} \Om)(e_i,X) = {\rm Ric}(JX,Z) -
{\rm Ric}^{*}(JX,Z).
\end{equation}
Symmetrising (\ref{ric-ric*-gen}) in $X$ and $Z$, one obtains the
relation between ${\rm Ric}''$ and the Nijenhuis tensor:
$$ (J {\rm Ric}'')_{X,Z} = \frac{1}{2} \sum_{i=1}^{4}
\Big( (\na^2_{e_i,Z} \Om)(e_i,X) + (\na^2_{e_i,X} \Om)(e_i,Z)
\Big) .$$ Anti-symmetrising (\ref{ric-ric*-gen}) in $X$ and $Z$
and using $d \Om = 0$, one gets
\begin{equation} \label{weitzOm}
\frac{1}{2}\na^*\na \Om = \rho_* -\rho ,
\end{equation}
which can be recognized as being nothing but the usual
Weitzenb\"ock formula \cite{besse} applied to the harmonic 2-form
$\Om$. From (\ref{weitzOm}) we also get immediately
$$ 2\rho_* '' = (\na^*\na \Om)'', \ \ \rho_0 = (\rho_* ')_0  \mbox
{ and } $$
\begin{equation} \label{s*-s}
 s^* - s = |\na \Om|^2, \mbox{ where } \; s^* = \frac{2\kappa +
 s}{3}
\end{equation}
is the trace of the star-Ricci tensor ${\rm Ric}^*$.

Finally, the component $W^+_3$ is determined by the $\leftr
\La^{0,2}M \rightr$-component of the relation (\ref{ricid}).

\vspace{0.2cm}

As $\na \Om$ is a section of the bundle $\leftr \La^{0,1}M
\otimes \La^{0,2}M  \rightr$, it will be very convenient at times
to express the higher jets of $\Om$ in terms of a (local) section
$\phi \in \leftr \La^{0,2}M \rightr$ of square-norm 2. Observe
that we have an $S^1$-freedom for the choice of $\phi$ and, for
this reason, such a form $\phi$ will be called a {\it gauge}.

Thus, assuming that we have made a choice for the gauge $\phi$, we
can write:
\begin{equation}\label{na-om}
\na\Om = a\otimes \phi - Ja \otimes J\phi.
\end{equation}
The gauge-dependent 1-form $a$ is dual to a vector field in
$D^{\perp}$ and satisfies $|\na \Om|^2 = 4|a|^2$, which shows that
$|a|^2$ is a gauge-independent quantity. As $\{
\frac{1}{\sqrt{2}}\Om, \frac{1}{\sqrt{2}}\phi ,
\frac{1}{\sqrt{2}}J\phi \}$ is an orthonormal basis of $\La^+M$,
the covariant derivatives of $\phi$ and $J\phi$ are given by
\begin{equation}\label{na-phi}
\na \phi = - a\otimes \Om + b\otimes J\phi; \ \ \na J\phi =
Ja\otimes \Om - b\otimes \phi,
\end{equation}
for some gauge-dependent 1-form $b$. In fact, if we change the
gauge  by
$$\phi' = (\cos\theta)\phi + (\sin\theta)J\phi,$$ the corresponding
1-forms change as follows
$$ a' = (\cos\theta)a - (\sin\theta)Ja;  \ b' = b + d\theta.$$
From (\ref{na-om}) and (\ref{na-phi}), we get
$$\na ^2|_{\Lambda^2 M}\Om =(da - Ja\wedge b)\otimes \phi -
(d(Ja) + a \wedge b)\otimes J\phi, $$ so, the Ricci identity for
$\Om$, (\ref{ricid}), can be rewritten as
\begin{equation}\label{*}
da - Ja\wedge b = - R(J\phi); \ d(Ja) + a\wedge b = - R(\phi).
\end{equation}
From the Ricci identity for $\phi$
$$(\na^2_{X,Y} - \na ^2_{Y,X})(\phi)  = [\phi,R_{X,Y}],$$
 after using (\ref{na-om}) and (\ref{na-phi}),
one gets the gauge independent relation
\begin{equation}\label{**}
db = a\wedge Ja -R(\Om) = a\wedge Ja - \rho_* .
\end{equation}
Note that  the closed, gauge independent 2-form
\begin{equation}\label{chernform}
\gamma_J = -db = R(\Om) - a\wedge Ja
\end{equation}
is a deRham representative of $2\pi
c_1(M,J)$, calculated with respect to the {\it first canonical
connection} of $(M,g,J)$
$$\na^1_X Y = \na_X Y - \frac{1}{2}J(\na_X J)(Y), $$
see e.g. \cite{gauduchon1}. We shall refer to $\gamma_J$ as the {\it canonical
Chern form} of $(M,g,J)$.

\subsection{The differential Bianchi identity for almost K\"ahler
4-manifolds} For a Riemannian 4-manifold $(M,g)$ the differential
Bianchi identity can be written as $ \delta W  = C $, where $W$
is the Weyl tensor, $C$ denotes the {\it Cotton-York} tensor of
$(M,g)$ and the co-differential $\delta$ acts on $W$ as on a
$\La^2 M$-valued 2-form. Recall that $C$ is defined by
$$C _{X, Y, Z} = - (d^{\na} h)(X,Y,Z) = - (\nabla _X h )(Y, Z) + (\nabla_Y
h)(X, Z),$$ where $h = \frac{1}{2}  {\rm Ric}_0 + \frac{s}{24} g
$ denotes the {\it normalized Ricci tensor} and $d^{\na}$ is the
Riemannian differential acting on $T^*M$ valued 1-forms. Both
$\delta W$ and $C$ are sections of the bundle $\La^1 M \otimes
\La^2 M$ and because of the splitting $\La^2 M = \La^+ M \oplus
\La^- M$, the differential Bianchi identity splits in two halves,
self-dual and anti-self-dual:
\begin{equation*}\label{bianchipm}
 \delta W ^+  = C ^+, \ \  \delta W ^-  = C ^- ,
\end{equation*}
where the $\pm$-superscript denotes the self-dual, resp.
anti-self-dual component of the corresponding tensor. We would
like to express the above formulae in terms of the various $U(2)$
curvature components of an almost K\"{a}hler 4-manifold $(M, g,
J,\Om)$. We start with the Cotton-York tensor $C$.

Contracting the differential Bianchi identity $\delta W = C$, we
obtain
$$\delta ({\rm Ric}_0 - \frac{s}{4} g) = 0$$
 and this,
together with the fact that the fundamental form $\Om$ is
harmonic, implies:
\begin{equation}\label{riccibianchi}
- \delta (\rho_0 - \frac{s}{4} \Omega) = J \delta ({\rm Ric}^{''}) .
\end{equation}
Taking the Hodge star operator of both sides, the above relation
is equivalent to
\begin{equation*} \label{drho}
 d \rho = * (J\delta {\rm Ric}^{''}) .
\end{equation*}
In particular, if ${\rm Ric}''=0$, then $\rho$ is closed \cite{draghici0}.

For a given vector field $Z$, we denote by $C_Z$  the section of
$\La^2 M$,  defined  by  $C_Z (X,Y): = C(X,Y,Z)$, and similarly we
define $C^+_Z$ and $C^-_Z$. Denote also by $A_Z$ the
$\La^2M$-valued 1-form given by
$$A_Z  = (d^{\nabla} {\rm Ric}^{''})_{Z} + \iota_{JZ}(d\rho),$$
where $\iota$ stands for interior derivative, and let
$A^{\pm}_Z$ be the self-dual and  anti-self-dual components of
$A_Z$.

\begin{Lemma} \label{lem3} Let $(M, g, J,
\Om)$  be an almost \ka 4-manifold. Then, for any vector field $Z$, the
Cotton-York tensor $C_Z$ is given by
\begin{equation}\label{cotton-york}
2C_Z = \nabla_{JZ} \rho_0 - \frac{1}{4}((d^Js)_Z) \Omega +
\nabla_{\rho_0(Z)} \Omega + \frac{1}{6} ds \wedge Z^{\flat} -
A_Z,
\end{equation}
where $Z^{\flat}$ denotes the $g$-dual 1-form of $Z$.
\end{Lemma}
\noindent {\it Proof.} Since $h = \frac{1}{2} {\rm Ric} -
\frac{s}{12} g$, the Cotton-York tensor  is written as:
\begin{equation}\label{basicpart}
C_Z = - \frac{1}{2} (d^{\na}{\rm Ric}^{''})_Z - \frac{1}{2} (d^{\na}
{\rm Ric}^{'})_Z + \frac{1}{12} ds\wedge Z^{\flat}.
\end{equation}
Taking into account that $\rho(\cdot, \cdot) = {\rm Ric}^{'}
(J\cdot, \cdot)$, we have
\begin{equation} \label{dricinv} \begin{split}
 (d^{\na} {\rm Ric}^{'})(X, Y, Z) &= (\nabla _X {\rm Ric}^{'})(Y, Z) -  (\nabla
_Y {\rm Ric}^{'}) (X, Z) \\ &= (d^{\na} \rho)(X,Y,JZ) + \Big(
\rho(Y, (\na_X J)Z) - \rho(X, (\na_Y J)Z) \Big) .
\end{split} \end{equation}
For the term $d^{\na} \rho$ we have
\begin{equation}\label{dricinvterm1} \begin{split}
(d^{\na} \rho)(X,Y,JZ) &= (\na_X \rho)(Y,JZ) - (\na_Y \rho)(X,JZ)
\\ &= - (\na_{JZ} \rho)(X,Y) + (d\rho)(X,Y,JZ) \\
&= - (\na_{JZ} \rho)(X,Y) + \iota_{JZ}(d\rho)(X,Y).
\end{split} \end{equation}

To refine the last term of (\ref{dricinv}), note that as an
algebraic object, $\na_X J$ is a skew-symmetric endomorphism of
$TM$, associated (by $g$-duality) to the section $\na_X \Omega$ of
the bundle $\leftr   \Lambda^{0,2} M \rightr$.  Thus
$\na_X J$ anti-commutes with $J$, and commutes with any
skew-symmetric endomorphism associated to a section of $\La^- M$;
in particular, it commutes with the endomorphism corresponding to
$\rho_0$ via the metric (which will be still denoted by
$\rho_0$). We  thus obtain
\begin{equation} \label{dricinvterm2} \begin{split}
\rho(Y, (\na_X J)(Z)) &- \rho(X, (\na_Y J)(Z) )  =
\frac{s}{4}\Big((\na_Y \Omega)(X, JZ) -(\na_X \Omega)(Y, JZ) \Big
) \\ & \ \ + (\na_X \Omega)(Y, \rho_0(Z)) - (\na_Y \Omega) (X,
\rho_0(Z))  \\ & \ \ = \frac{s}{4}(\na_{JZ} \Omega) (X,Y) -
(\na_{{\rho_0} (Z)} \Omega)(X,Y),
\end{split} \end{equation}
where for the last equality we used the closedness of $\Omega$.

Substituting (\ref{dricinvterm1}) and (\ref{dricinvterm2}) in
(\ref{dricinv}), we get
\begin{equation} \label{dricinv+} \begin{split}
(d^{\na} {\rm Ric}^{'})_Z &= - \na_{JZ} \rho_0 - \na_{\rho_0(Z)}
\Om \\ & \ \ \ + \frac{1}{4}((d^Js)_Z) \Omega + \iota_{JZ}(d\rho).
\end{split} \end{equation}
From (\ref{basicpart}) and (\ref{dricinv+}), we obtain relation
(\ref{cotton-york}) claimed in the statement. $\Box$

\begin{Lemma} \label{Bianchi-lem} Let $(M, g, J, \Om)$ be an
almost K{\"a}hler 4-manifold. Then the differential Bianchi
identities $\delta W ^{+}  = C ^{+}, \ \delta W ^{-}  = C ^{-} $
are equivalent to:
\begin{eqnarray}\label{bianchi+}
 0 &=& -\frac{1}{4} (d^J(\kappa - s))_Z \Omega + \frac{1}{6} (d(\kappa
- s) \wedge Z^{\flat})^+ + (\delta \rho_*'')_Z \Omega  \\
\nonumber & & + \frac{\kappa}{4} \nabla_{JZ} \Omega -
\nabla_{\rho_0(Z)} \Omega + \nabla_{\rho_{*}''(Z)} \Omega +
\nabla_{JZ} \rho_{*}'' + 2(\delta W_3^+)_Z + A^+_Z  \; ; \\
\label{bianchi-} 0 &=& \nabla_{JZ} \rho_0 + \frac{1}{6} (d s
\wedge Z^{\flat})^- - 2\delta W^-_{Z} - A^-_Z \; .
\end{eqnarray}
\end{Lemma}
\noindent{\it Proof.} Relation (\ref{bianchi-}) follows from $\delta W
^{-} = C ^{-}$ by simply taking the anti-self-dual component of
(\ref{cotton-york}). Similarly, for (\ref{bianchi+}) we take the
self-dual component of (\ref{cotton-york}), but we also use the
relations
\begin{equation}\label{deltaw^+1}
(\delta W_1^+)_Z = -\frac{1}{8} (d^J\kappa)_Z \Om + \frac{\kappa}{8}
\na_{JZ} \Om + \frac{1}{12} (d\kappa \wedge Z^{\flat})^+ ,
\end{equation}
\begin{equation}\label{deltaw^+2}
(\delta W^+_2)_Z = \frac{1}{2} \nabla_{JZ} (\rho_*)''  +
\frac{1}{2} \nabla_{(\rho_*)''(Z)} \Omega + \frac{1}{2} (\delta
\rho_*'')(Z) \Omega ,
\end{equation}
which  easily follow from (\ref{w^+1}) and
(\ref{w^+2}).
$\Box$

\vspace{0.2cm}

It will be useful to further determine the $\Om$-component and the
$\leftr
\Lambda^{0,2} M \rightr$-component of (\ref{bianchi+}) in accordance
with the decomposition (\ref{U(2)La+}) of $\La^+ M$. These are given by
\begin{cor} \label{bianchi+rem} For any almost \ka 4-manifold the
following identities hold:
\begin{eqnarray}\label{bianchiOm}
 0 &=& -\frac{1}{3} (d^J(\kappa - s))_Z + 2\langle (\delta W_3^+)_Z,
\Omega \rangle
\\ \nonumber
 & & + 2 (\delta \rho_*'')_Z -  \langle \rho_*'',
\nabla_{JZ} \Omega \rangle + \langle A^+_Z, \Omega \rangle  \; ;
\\
\label{bianchi+2,0} 0 &=& \frac{1}{6} (d(\kappa - s) \wedge
Z^{\flat})^{''} + \frac{\kappa}{4} \nabla_{JZ} \Omega -
\nabla_{\rho_0(Z)} \Omega + 2(\delta W_3^+)_Z^{''} \\ \nonumber
& &+ \nabla_{\rho_*''(Z)} \Omega + (\nabla_{JZ}\rho_*
'')^{''}  + (A^+_Z)^{''}.
\end{eqnarray}
\end{cor}

\subsection{A Weitzenb\"ock formula and unique continuation
of the Nijenhuis tensor} Recall that the {\it weak unique
continuation property} for a map $u : M \to E$ between connected
Riemannian manifolds, $M$ and $E$, says that if $u$ is a constant
map on an open subset $U \subset M$, then $u$ is constant
everywhere; for the {\it strong} unique continuation property,
the condition that $u$ is constant on an open subset is replaced
with the assumption that at a given point, $u$ has a contact of
infinite order with the constant map (see e.g. \cite{kazdan}).

In this subsection we shall show that the (strong) unique
continuation property holds for the Nijenhuis tensor of an almost
\ka 4-manifold which satisfies certain curvature conditions. In
fact, we shall concentrate our attention on $\na \Om$, which is
identified to the Nijenhuis tensor via (\ref{naOm-N}). As
typically happens in Riemannian geometry, the unique continuation
property appears as a consequence of a Weitzenb\"ock formula. In
our case, we shall use a general Weitzenb\"ock formula of
Bourguignon \cite{bourg} applied to $\na \Om$, seen as a section
of the bundle of $\La^2 M$-valued 1-forms. The role of the
``constant map'' in this setting is played by the zero section.

\begin{prop}\label{ucp} Let $(M,g,J,\Om)$ be an almost \ka
4-manifold whose curvature satisfies the third Gray condition {\rm
(\ref{gray3})}. Then, the  following identity holds:
\begin{equation}\label{weitzenbock1}
\na^* \na (\na \Om) + \frac{3s}{4} \na \Om + \na_{{\rm
Ric}_0(\cdot)} \Om =0,
\end{equation}
In particular, if $M$ is connected, then  the (strong) unique
continuation property holds for $\na \Om$ (hence for the
Nijenhuis tensor $N$ as well).

\end{prop}

\noindent {\it Proof:} We shall in fact establish a formula for
$\na^* \na (\na \Om)$ on an arbitrary almost \ka 4-manifold, and
then relation (\ref{weitzenbock1}) will be an immediate
consequence. The starting point is the general Weitzenb\"ock
formula for a $\La^2 M $-valued 1-form $V$ (e.g. see
\cite{bourg}, p.282):
\begin{equation}\label{jpb2}
\Big( (d^{\nabla}\delta^{\nabla} +  \delta^{\nabla}d^{\na}) V
\Big)_X = \Big( \na^*\na V \Big)_X + V_{{\rm Ric}(X)} + (R \cdot
V)_X .
\end{equation}
In the above formula
\begin{enumerate}
\item[$\bullet$] $d^{\na}$ is the Riemannian differential acting on
$\La^2 M$-valued $k$-forms;

\item[$\bullet$] $\delta^{\na}$ is the formal adjoint of
$d^{\na}$;

\item[$\bullet$] the curvature action is given by:
\begin{equation*}
(R \cdot V)_X = \sum_{i=1}^4 [R_{X, e_i} , V_{e_i}],
\end{equation*}
\end{enumerate}
where $\{e_1, ... , e_4 \}$ is an orthonormal basis of $TM$, the
2-forms $R_{X,e_i}$, $V_{e_i}$ are freely identified with the
corresponding skew-symmetric endomorphism of $TM$ and
$[\cdot,\cdot]$ denotes the commutator.
We apply (\ref{jpb2}) to $V = \na \Om$; using (\ref{weitzOm}),
(\ref{s*-s}) and (\ref{ricid}) we get
\begin{eqnarray} \label{delnaom}
\delta^{\na} (\na \Om) &=& \na^* \na \Om = 2(\rho_* - \rho) = 2
\rho_*'' + \frac{(\kappa -s)}{3} \Om \; ,
\\ \label{dnaom}
(d^{\na} (\na \Om))_{X,Y} &=& (\na^2_{X,Y} \Om) - (\na^2_{Y,X}
\Om) = [J, R_{X,Y}]. \;
\end{eqnarray}
From (\ref{delnaom}) we derive immediately
\begin{equation} \label{ddelnaom}
(d^{\na} \delta^{\na} (\na \Om))_X = 2 \na_X \rho_*'' +
\frac{1}{3} d(\kappa - s)_X \Om + \frac{(\kappa-s)}{3} \na_X \Om \; ,
\end{equation}
whereas a short computation starting from (\ref{dnaom}) leads to
\begin{equation} \label{deldnaom1}
 (\delta^{\na} d^{\na} (\na \Om))_X = - (R \cdot (\na \Om))_X +
 \Big( (\na_{e_i} R)_{e_i,X}(J\cdot, \cdot) +
 (\na_{e_i} R)_{e_i,X}(\cdot, J\cdot) \Big) \; .
\end{equation}
Using the differential Bianchi identity and (\ref{dricinv+}), we
get
\begin{equation} \begin{split} \label{t2deldnaom}
\Big( (\na_{e_i} R)_{e_i,X}(J\cdot, \cdot) &+
 (\na_{e_i} R)_{e_i,X}(\cdot, J\cdot) \Big) = (d^{\na} {\rm
 Ric})_X (J\cdot, \cdot) + (d^{\na} {\rm Ric})_X (\cdot, J\cdot) \\
 &= A_X (J\cdot, \cdot) + A_X (\cdot, J\cdot)+2
(\na_{{\rm Ric}^{'}_0(X)} \Om)(\cdot,\cdot) \\
& = - 2 (J (A^+_X)^{''})(\cdot,\cdot) +2 (\na_{{\rm
Ric}^{'}_0(X)} \Om)(\cdot,\cdot) \;.
\end{split} \end{equation}
From the relations (\ref{ddelnaom}-\ref{t2deldnaom})
 and (\ref{jpb2}), we obtain
\begin{equation} \label{laplacenaom*} \begin{split}
\na^* \na (\na \Om) &= 2 \na \rho_*''- 2J (A^+)^{''}
- \na_{{\rm Ric}(\cdot)} \Om\\
& \ \ - 2R \cdot (\na \Om) +2 \na_{{\rm Ric}^{'}_0(\cdot)} \Om +
\frac{1}{3} d(\kappa - s) \otimes \Om + \frac{(\kappa-s)}{2} \na \Om \;.
\end{split} \end{equation}
 It remains to detail the expression of the term $R \cdot (\na
\Om)$. For this, we shall use the $U(2)$-decomposition of the
curvature (\ref{u(2)}) and compute the action of each component
on $\na \Om$. It is clear that the action of $W^-$ is zero, and
for the action of the component $\widetilde{{\rm Ric}_0^{''}}$ we
will not attempt to do any simplification, as ${\rm Ric}_0''=0$
for almost \ka 4-manifolds  satisfying (\ref{gray3}). For the
other components though, from their definitions (see
(\ref{w^+1}), (\ref{w^+2}), (\ref{w^+3})),  one computes
successively:
\begin{eqnarray*}
\Big( \frac{s}{12} {\rm Id}_{| \La^2 M} \cdot (\na \Om) \Big)_X &=&
\frac{s}{12} (\na_X \Om) ;\\
\Big( W_1^+ \cdot (\na \Om) \Big)_X &=&
\frac{\kappa}{6} (\na_X \Om) ;\\
\Big( W_2^+ \cdot (\na \Om) \Big)_X &=& - (\na_{J\rho_*''(X)}
\Om) - \frac{1}{2} \langle \rho_*'', \na_X \Om \rangle \Om ;\\
\Big( W_3^+ \cdot (\na \Om) \Big)_X &=& \langle (\delta W_3^+)_{JX},
\Om \rangle \Om ;\\
\Big( \widetilde{{\rm Ric}_0^{'}} \cdot (\na \Om) \Big)_X &=&
\na_{{\rm Ric}_0^{'}(X)} \Om \; .
\end{eqnarray*}
Finally, using the relations above, back in (\ref{laplacenaom*}),
we obtain:
\begin{equation} \label{laplacenaom} \begin{split}
\Big( \na^* \na (\na \Om) \Big)_X &= - \na_{{\rm Ric}_0^{'}(X)}
\Om - \frac{3s}{4} \na_X \Om \\
& \ \ \ + \Big( \frac{1}{3} d(\kappa - s)_X - 2\langle (\delta
W_3^+)_{JX}, \Om \rangle \Big) \Om \\
& \ \ \ + 2(\na_X \rho_* '')'' + 2
\na_{J\rho_*''(X)} \Om \\
& \ \ \ - \na_{{\rm Ric}_0^{''}(X)} \Om - 2 J(A_X)^{''} -
2\Big( \widetilde{{\rm Ric}_0^{''}} \cdot (\na \Om) \Big)_X \; .
\end{split} \end{equation}

For an almost \ka 4-manifold which satisfies the third Gray
condition, the terms in the last two lines of the above equation
all vanish, as they contain ${\rm Ric}^{''}$ and $\rho_* ''$, see Lemma
\ref{lem1}. In fact, the expression of the second line is also $0$, as a
consequence of the Bianchi relation (\ref{bianchiOm}). Thus,
under the assumption ${\rm Ric}^{''}=0 \; \& \; \rho_* ''=0$, the
general formula (\ref{laplacenaom}) reduces to
(\ref{weitzenbock1}); the fact that the relation (\ref{weitzenbock1})
implies the (strong) unique continuation property for $\na \Om$ (and
therefore also for the Nijenhuis tensor $N$), follows from the classical
result of Aronszajn
\cite{arons}. $\Box$

\vspace{0.2cm}
\begin{rem}\label{weitzenbock-rem}
{\rm (i) Note that solely for the purpose of getting the unique
continuation property for $\na \Om$, we need not make the effort
to compute the action of each of the curvature components on $\na
\Om$. Indeed, to apply the result of Aronszajn \cite{arons}, it is
enough to obtain an estimate of the form
$$ |\na^* \na (\na \Om) |^2 \leq M \Big( |\na (\na \Om)|^2 + |\na \Om|^2 \Big)
\; , $$ and relation (\ref{laplacenaom*}) already implies such an
estimate, once we also note that
$$ \frac{2}{3}|d(\kappa - s)| =
|d(|\na \Om|^2)| \leq 2|\na (\na \Om)| |\na \Om| =
\frac{4}{3}(\kappa - s)|\na (\na \Om)|.$$ More generally, from
(\ref{laplacenaom}) (or equally well from (\ref{laplacenaom*})),
one can deduce that the (strong) unique continuation property for
$\na \Om$ holds for almost \ka 4-manifolds for which two of the
three tensors, ${\rm Ric}^{''}, \rho_* ''$, and $W_3^+$ vanish.
Indeed, in the case when  $\rho_* '' = 0  \ \&  \ W^+_3 = 0$, the
differential Bianchi identity (\ref{bianchi+2,0}) can be written
as
$$ - (A^+_Z)^{''} = \frac{1}{6} (d(\kappa - s) \wedge
Z^{\flat})^{''} + \frac{\kappa}{4} \nabla_{JZ} \Omega -
\nabla_{\rho_0(Z)} \Omega . $$ Putting the above relation back in
(\ref{laplacenaom*}), one again obtains the needed estimate to
apply \cite{arons}.  One can proceed similarly in the case ${\rm
Ric}^{''} = 0 \ \& \ W^+_3=0$. In fact, in this latter case,
employing a different estimate, it was
shown in \cite{AA} that an even stronger version of the unique
continuation property of $\na \Om$ holds; namely, if $\na \Om$
vanishes at one point, then it must be identically zero.

(ii) One would expect formula (\ref{laplacenaom}) to eventually
lead to the integrability of compact almost K\"ahler 4-manifolds
satisfying certain curvature conditions. This is indeed the case:
along these lines it was proved in \cite{Dr2} that on a compact
symplectic 4-manifold there are no critical, strictly almost \ka
structures of everywhere non-negative scalar curvature.}
\end{rem}

\section{The examples of almost \ka 4-manifolds satisfying the
third Gray condition}

We give here more details about the explicit construction given
in Theorem \ref{thm1}; the material in this section has appeared
in \cite{AAD}.

Let $(\Sigma, g_{\Sigma})$ be any oriented Riemann surface and $h=w+iv$
be  a holomorphic function
on $\Sigma$, such that $w= {\frak Re}(h)$ is a positive (harmonic)
function on $\Sigma$.
We can then locally write the metric $g_{\Sigma}$ as
$$g_{\Sigma}=e^{u}w(dx^{\otimes 2} + dy^{\otimes 2}),$$ where $x,y$ are
(local) isotherm coordinates of $(\Sigma, g_{\Sigma})$ and $u$ is a
function on $\Sigma$. We consider the metric
\begin{eqnarray}\label{g}
g &=& g_{\Sigma} + w dz^{\otimes 2} + \frac{1}{w}(dt + vdz)^{\otimes 2}\\
\nonumber
&=& e^{u}w(dx^{\otimes 2} + dy^{\otimes 2})+ w dz^{\otimes 2} +
\frac{1}{w}(dt + vdz)^{\otimes 2},
\end{eqnarray}
defined  on $M= \Sigma \times {\mathbb R}^2$, where $(z,t)$ are
the canonical coordinates of ${\mathbb R}^2$. Clearly, $K_1=
\frac{\partial}{\partial z}$ and $K_2 = \frac{\partial}{\partial
t}$ are two commuting Killing fields, so that $g$ is an ${\mathbb
R}^2$-invariant metric compatible with the product structure in
the sense of \cite{joyce}.

As
observed in
\cite{LeB}, the metric (\ref{g}) carries a compatible \ka structure, $I$,
whose fundamental form is
$${\overline \Om} =  \Omega_{\Sigma} + dz\wedge (dt + vdz) = (e^uw)
dx\wedge dy + dz\wedge dt,$$ where $\Omega_{\Sigma}$ denotes the volume
form of
$(\Sigma, g_{\Sigma})$. Equivalently, $I$ is defined by its action on
1-forms:
\begin{equation}\label{I}
I(dx)= dy; \ I(dz) = \frac{1}{w}(dt + vdz).
\end{equation}
 Thus, $I$
is compatible with the product structure as well, and induces an
orientation on  $M$ and on each of the factors $\Sigma$ and
${\mathbb R}^2$. Besides $I$, one can consider the almost complex
structure $J$, which coincides with $I$ on $\Sigma$, but which is
equal to $-I$ on ${\mathbb R}^2$, i.e.
\begin{equation}\label{J}
J (dx) = dy;  \ J(dz) = - \frac{1}{w}(dt + vdz).
\end{equation}
Thus, $J$ is compatible with $g$ and yields on $M$ the
orientation opposite to the one induced by $I$. Furthermore, the
fundamental form of $J$ is given by
$$\Omega = \Omega_{\Sigma} - dz\wedge (dt + vdz) = \Omega_{\Sigma} -
dz\wedge dt,$$ and is clearly closed, meaning that $(g, J)$ is an almost
\ka structure.  It is easily seen that $J$ is integrable (i.e.
$(g,J)$ is \ka) if and only if $h$ is  constant (i.e. $g$ is a
product metric), a possibility which we exclude in what follows below.

An important feature of the construction comes from the following
observation: at any point where $dw \neq 0$, the \ka nullity ${D}
= \{ TM \ni X : \na_X J = 0 \}$ of $J$ is a two dimensional
subspace of $TM$, which is tangent to the surface $\Sigma$ (see
\cite{AAD}); thus, the \ka structure $I$ and the almost \ka
structure $J$ are equivalently related by
$$J_{D} = I_{D} \ \ J_{{D}^{\perp}} = - I_{{D}^{\perp}}.$$

The almost \ka structure $(g,J, \Om)$ is sufficiently explicit
to make it straightforward, though tedious, to check that its
curvature satisfies the third Gray condition (\ref{gray3}); see
our previous paper \cite{AAD} for detailed calculations. But the
heart of \cite{AAD} consists of the following abstract
characterization of the almost \ka structures given by
(\ref{g}-\ref{J}):
\begin{theo}\label{thm3}{\rm \cite[Theorem 2]{AAD}} Let $(M,g,J,\Om)$ be a
strictly almost \ka 4-manifold whose curvature satisfies
{\rm (\ref{gray3})}. At each point $x \in M$ where the Nijenhuis
tensor of $J$ does not vanish, consider the orthogonal almost
complex structure $I$ which is equal to $J$ on the \ka nullity $D_x
\subset T_xM$, but to $-J$ on the orthogonal complement of ${D}_x$.
If we suppose that $(g,I)$ is \ka on an open subset $U$, then for every
sufficiently small neighborhood of $x \in U$,  $(g,J,I)$  is given by
{\rm (\ref{g}-\ref{J})}, up to a homothety.
\end{theo}

The key point for obtaining the above result in \cite{AAD} was to
show that the distribution $D^{\perp}$ is (locally) spanned by
commuting hamiltonian Killing fields, a subject to a
straightforward verification of the integrability condition of a
relevant Frobenius system; then, according  to \cite{LeB}, the
\ka metric $(g,I)$ can be put (locally) in the form
(\ref{g}-\ref{I}) and Theorem \ref{thm3}  follows easily.

On the other hand, it was shown in \cite{Arm1} that  when $(g,J)$
is an Einstein, strictly almost \ka structure satisfying
(\ref{gray3}), $(g,I)$ is necessarily \ka; the same was later
established for strictly almost \ka 4-manifolds satisfying the
second curvature condition of Gray (see \cite{AAD}). In the next
section we generalize these results.

\section{Proof of Theorem 1}

The following
proposition, which via Proposition \ref{ucp} and Theorem \ref{thm3}
implies Theorem
\ref{thm1}, is the technical core of the present paper.
\begin{prop}\label{localadd-prop}
Let $(M,g,J,\Om)$ be a strictly almost K{\"a}hler 4-manifold
satisfying the third Gray condition {\rm (\ref{gray3})} and let
$(g,I)$ be the opposite almost Hermitian structure associated to
$(g,J)$ as in Theorem \ref{thm3}. Then, on the open set $U
\subset M $, where $|\na J| \neq 0$, $(g,I)$ is
a K{\"a}hler structure.
\end{prop}
There are three main steps in the proof, corresponding roughly to
Lemmas \ref{step1}, \ref{step2} and  \ref{step3} below. In Lemma
\ref{step1}, we first establish that under the third Gray
condition (\ref{gray3}) for $(g,J)$, the opposite almost Hermitian
structure $(g,I)$ is {\it almost \ka} and we specifically find, in
accordance with Lemma \ref{lem2}, the exact form representing the
difference $\gamma_{I} - 3\gamma_{J}$, where $\gamma_J$ and
$\gamma_I$ are the canonical Chern forms of $(g,J)$ and $(g,I)$
respectively (see Section 2.3). Under the {\it second} Gray
condition studied in \cite{AAD}, or with the Einstein assumption
made in \cite{Arm1}, Lemma \ref{step1} is essentially enough to
obtain that  $(g,I)$ is K\"ahler.  Assuming the (weaker) third
Gray condition only, the same conclusion is still true, but not so
immediate.  The idea of the proof is to proceed by contradiction,
thus assuming  that $(g,I)$ is {\it not} \ka and to try to iterate
the construction. This is somehow complicated by the fact that,
{\it a priori}, $(g,I)$ does not satisfy the third Gray condition
(\ref{gray3}). Nevertheless, in Lemma \ref{step2} we show that the
curvature has the necessary features to allow the iteration of the
construction. The last step is essentially done in Lemma
\ref{step3}: We consider $(g, {\tilde J})$, the opposite almost
Hermitian structure associated to $(g,I)$ and, analogously to
Lemma \ref{step1}, we  get an explicit expression for the exact
form representing $\gamma_{{\tilde J}} - 3\gamma_I$, where
$\gamma_{\tilde J}$ denotes the canonical Chern form of
$(g,{\tilde J})$. But one immediately notices that, in view of
Lemma \ref{step1},  the ``opposite of the opposite'' almost
complex structure ${\tilde J}$ is nothing else than $-J$, hence we
also have $\gamma_{{\tilde J}} = - \gamma_J$. Therefore, the
expression for $\gamma_{{\tilde J}} - 3\gamma_I$ obtained in Lemma
\ref{step3} provides another relation between $\gamma_J$ and
$\gamma_I$, {\it different} than the one obtained in Lemma
\ref{step1} (see also Lemma \ref{lem2}). The combination of the
two leads to the desired contradiction.

This is the outline of the proof of Proposition
\ref{localadd-prop}. Before we proceed to detail each of the
three mentioned steps, let us first recall some notations from
Section 2 and introduce some new ones, for convenience. As
mentioned in Section 2.3, it is  useful to make computations
using a gauge $\phi$; recall that $\phi$ is a local section of
norm $\sqrt{2}$ of the bundle $\leftr \La^{0,2}M \rightr $. We
thus choose an arbitrary gauge $\phi$ and fix it. By formulae
(\ref{na-om}) and (\ref{na-phi}), we define the local 1-forms $a$
and $b$ corresponding to $\phi$. We shall denote the (local)
endomorphisms of $TM$ induced by $\phi$ and $J \phi$ by $J_1$
and, respectively, $J_2$; note that the triple $(J, J_1, J_2 = J
J_1)$ is a local almost quaternionic structure on the manifold.
We also make the convention to denote a tensor (or form) defined
by the opposite almost complex structure $(g, I)$ with the same
symbol as the corresponding one defined by $(g, J)$, but with a
``bar'' sign on top. Thus, for instance, the fundamental form of
$(g,I)$ is denoted by ${\overline \Om}$, the Ricci form of
$(g,I)$ by ${\bar \rho}$, etc.

Because of the existence of the opposite  almost Hermitian structure
$(g,I)$, the decomposition (\ref{Lambda2r}) of the bundle of 2-forms
further refines to:
\begin{equation*}\label{Lambda2rIJ}
\La^2M = {\mathbb R}\cdot \Om \oplus \leftr \La_J^{0,2} M \rightr
\oplus {\mathbb R}\cdot {\overline \Om} \oplus \leftr
\La_I^{0,2} M \rightr ,
\end{equation*}
where $\leftr \La_J^{0,2} M \rightr$, $\leftr \La_I^{0,2} M
\rightr$ denote the underlying real bundles of the anti-canonical
bundles of $(g,J)$, $(g,I)$, respectively. We shall keep the
convention made earlier to use the superscripts $'$ and $''$ for
the $J$-invariant, respectively, $J$-anti-invariant parts of a
2-form, or symmetric 2-tensor. But we shall need occasionally to
also consider $\leftr \La_I^{0,2} M \rightr$-components of
2-forms and we make the convention to use the notation $( \cdot
)''_I$ for this. Thus, for instance, $({\bar \rho}_{*})''_I$ will
denote the $\leftr \La_I^{0,2} M\rightr$-component of the
$*$-Ricci form ${\bar \rho}_{*}$ of $(g,I)$. In addition to
$d^J$, we shall also use $d^I$, the complex differential with
respect to $I$, acting on functions  by $d^I f = I df$ (see Sec.
2.3).

\vspace{0.2cm}

The Lemma below was essentially proved in our previous work
\cite{AAD}.

\begin{Lemma}\label{step1} Let $(M,g,J,\Om)$ be a strictly almost \ka
4-manifold satisfying {\rm (\ref{gray3})}. Then, on the open set $U$
where $N \neq 0$, the \ka nullity $D$ of $(g,J)$, and its
orthogonal complement, $D^{\perp}$, are both involutive
distributions. Furthermore, on $U$,
$(g,I)$ is an almost \ka structure and the \ka nullity
of $(g,I)$ contains the distribution $D^{\perp}$.
Moreover, the following relation between the canonical Chern forms
$\gamma_I$ and $\gamma_J$ of $(g,I)$ and $(g,J)$
holds:
\begin{equation} \label{gammaI-gammaJ}
\gamma_I = 3\gamma_J - d d^J(\ln{|\na \Om|^2}).
\end{equation}
\end{Lemma}

\noindent
{\it Proof.} The opposite almost complex structure $I$ is defined (on $U$)
through the distributions $D$ and $D^{\perp}$, which in turn are
determined by the 1-jet of $J$. To obtain information about the
1-jet of $I$, we thus need to study the 2-jet of $J$
(equivalently, the 2-jet of $\Om$).
For this purpose, let us define the 1-forms $m_i,
n_i, \ i=1,2$, by the first relation below, while the other three
relations follow easily from the first one and (\ref{na-om}),
(\ref{na-phi}).
\begin{equation} \label{na-a} \begin{split}
\na a &= m_1\otimes a + n_1\otimes Ja + m_2\otimes J_1 a +
n_2\otimes J_2 a ; \\
\na(Ja) &= -n_1 \otimes a + m_1 \otimes Ja + (a-n_2)\otimes J_1 a
 + (m_2 - Ja) \otimes J_2 a; \\
\na (J_1 a) &= -m_2 \otimes a + (n_2 - a) \otimes Ja + m_1 \otimes
J_1 a  + (b - n_1) \otimes J_2 a; \\
\na (J_2 a) &= -n_2 \otimes a + (Ja - m_2) \otimes Ja + (n_1-b)
\otimes J_1 a + m_1 \otimes J_2 a. \end{split}
\end{equation}
Note that for a strictly almost \ka 4-manifold $(M, g, J, \Om)$
satisfying (\ref{gray3}), the identity (\ref{bianchiOm}) becomes
\begin{equation*} \label{bianchi1ak3}
 0 = \frac{1}{6} (d^J(\kappa - s))_Z - \langle (\delta W_3^+)_Z,
 \Omega \rangle.
\end{equation*}
Using (\ref{s*-s}), (\ref{w^+3}) and (\ref{na-om}),
this can be further written as
\begin{equation} \label{lem3-1}
\frac{1}{4} d |\na \Om|^2 = d |a|^2 = -J \langle \delta W_3^+,
\Om \rangle = - 2 J W^+_3(\phi)(a) .
\end{equation}
Hence the 1-form $m_1$ is immediately determined to be:
\begin{equation} \label{m1}
m_1 = \frac{1}{2}d(\ln{|\na \Om|^2}) = - \frac{1}{|a|^2}
JW_3^+(\phi)(a).
\end{equation}
As a consequence, observe that $m_1^{\sharp} \in D$, where here
and henceforth the superscript $\sharp$ denotes the dual vector
field of a 1-form, through the metric. The forms $n_1, m_2, n_2$
are determined by the relations (\ref{*}). Indeed, under the
condition (\ref{gray3}),  the 2-forms $R(\phi)$ and $R(J\phi)$
are $J$-anti-invariant, thus,  according to  (\ref{*}), we have
$$ (da - Ja \wedge b)' = 0 , \; \; \; (d(Ja) + a \wedge b)' = 0 .$$
These relations combined with (\ref{na-a}) lead to
\begin{equation} \label{naa} \begin{split}
n_1 &= -b - Jm_1 = -b - \frac{1}{2}d^J (\ln{|\na \Om|^2}) ; \\
m_2 &= \frac{1}{2}Ja + Jm_0 ; \\
n_2 &= -Jm_2 = \frac{1}{2}a + m_0, \end{split}
\end{equation}
where $m_0$ is a 1-form dual to a vector field in $D$.

\vspace{0.2cm}

Thus, formulae (\ref{na-a}), (\ref{m1}) and (\ref{naa}) imply
\begin{equation*} \begin{split}
da(X,X') & = d(Ja)(X,X') = 0, \; \forall X,X' \in D , \\
d(J_1 a)(Y,Y') & = d(J_2  a)(Y,Y') = 0 , \; \forall Y,Y' \in D^{\perp}.
\end{split} \end{equation*}
As the vector fields dual to $a, Ja$ generate $D^{\perp}$
and the vector fields dual to $J_1a, J_2a$ generate $D$, the above
relations show that $D^{\perp}$ and $D$ are involutive distributions.

\vspace{0.2cm}

To show that $(g,I)$ is almost K\"ahler, observe first that the
  fundamental forms of
$(g,I)$ and $(g,J)$ are related by ${\overline \Om} = \Om -
\frac{2}{|a|^2} a \wedge Ja $. We take the covariant derivative
of this relation and
 use (\ref{na-om}) and (\ref{na-a}-\ref{naa})
to get
\begin{equation} \label{na-bar-Om}
\na {\overline \Om} = 2m_0 \otimes {\bar \phi} - 2{I}m_0 \otimes
{I}{\bar \phi} ,
\end{equation}
where ${\bar \phi}  \in \leftr  \Lambda_I^{0,2} M \rightr$
is the natural gauge for $(g,I)$, determined by $\phi$ through
the relation
\begin{equation} \label{bargauge}
{\bar \phi} = \phi + \frac{2}{|a|^2} Ja \wedge J_2 a.
\end{equation}
Then, $d{\overline \Om} = 0$ is immediate from (\ref{na-bar-Om}),
and so is the claim about the \ka nullity of $(g,I)$ since   ${\bar a} =
2 m_0$ is dual to a vector field in $D$.

\vspace{0.2cm}

Finally, taking the covariant derivative of (\ref{bargauge}) and using
(\ref{na-phi}) and (\ref{na-a}-\ref{naa}), one finds
\begin{equation} \label{na-bar-phi}
 \na {\bar \phi} = {\bar b} \otimes I{\bar \phi}
- {\bar a} \otimes {\overline \Om} ,
\end{equation}
where ${\bar a} = 2m_0$  and
\begin{equation} \label{bar-b}
{\bar b} = 3b + d^J (\ln{|\na \Om|^2}) .
\end{equation}
Since $\gamma_I = - d{\bar b}$ and $\gamma_J = - d b$ (see
(\ref{chernform})), relation (\ref{gammaI-gammaJ}) follows. $\Box$

\vspace{.2cm}

The next step is to gather more information on the curvature
components of the almost \ka metric $(g,I)$.

\begin{Lemma}\label{step2} Let $(M,g,J)$ be a
strictly almost \ka 4-manifold satisfying {\rm (\ref{gray3})} and
$(g, I)$ be the opposite almost \ka structure as above. Then, the
trace-free Ricci form $\rho_0$ of $(g,J)$ and the $\leftr
\La_I^{2,0} M\rightr$-component of the star-Ricci form ${\bar
\rho}_{*}$ of $(g, I)$ are given by:
\begin{equation} \label{barpsi}
 \rho_0= \frac{1}{4}(s + \frac{ |\na {\overline \Om}|^2}{4}) {\overline \Om}
 + {\overline \Psi},
 \; \; \;  ({\bar \rho}_{*})''_I = {\overline \Psi}, \;
\mbox{ where } {\overline \Psi} = \frac{1}{2} \na_{d(\ln{|\na
\Om|^2})} {\overline \Om} .
\end{equation}
Moreover, the identity {\rm (\ref{bianchi-})} reads as
\begin{equation} \label{barlem3-1}
\frac{1}{4} d |\na {\overline \Om}|^2 = d |{\bar a}|^2 = - I
\langle \delta W_3^{-}, {\overline \Om} \rangle = - 2I
W_3^{-}({\bar \phi})({\bar a}) \; .
\end{equation}
\end{Lemma}

\noindent {\it Proof.} We first prove the equalities in
(\ref{barpsi}). With respect to the gauge ${\bar \phi}$ from the
previous lemma, the Ricci identity for ${\overline \Om}$ takes the
form (compare with (\ref{*}))
\begin{equation}\label{bar*}
d{\bar a} - I{\bar a}\wedge {\bar b} = - R(I{\bar \phi}); \
d(I{\bar a}) + {\bar a}\wedge {\bar b} = - R({\bar \phi}).
\end{equation}
We proved in Lemma \ref{step1} that ${\bar a}^{\sharp} \in D$ and
that $D$ is an involutive distribution. It thus follows that the
left-hand sides of the relations (\ref{bar*}) vanish on $X \wedge
JX$ with $X \in D^{\perp}$. This is equivalent with the fact that
the $\leftr \Lambda_I^{0,2} M\rightr$-components of $R(\Om)$ and
$R({\overline \Om})$ are equal. Thus $$ (\rho_0)''_I = ({\bar
\rho}_{*})''_I = {\overline \Psi} ,$$ and it remains to obtain the
claimed expressions of ${\overline \Psi}$ and the ${\overline
\Om}$-component of $\rho_0$. For this we turn to relation
(\ref{gammaI-gammaJ}). Using (\ref{**}) and its analogue for
${\bar \phi}$
\begin{equation*} \label{bar**}
d{\bar b} = {\bar a} \wedge I{\bar a} - R({\overline \Om}) ,
\end{equation*}
relation (\ref{gammaI-gammaJ}) is equivalent to
\begin{equation} \label{db-comparison}
\frac{3}{2} R(\Om) - \frac{1}{2} R({\overline \Om}) - \frac{1}{2}
dd^J(\ln{|\na \Om|^2}) - \frac{3}{2} a \wedge Ja + \frac{1}{2}
{\bar a} \wedge J{\bar a} = 0.
\end{equation}
Taking the $\leftr \Lambda_I^{0,2} M \rightr$-component of
(\ref{db-comparison}) we obtain
$$ {\overline \Psi} = \frac{1}{2} (dd^I (\ln{|\na \Om|^2}))''_I
 = \frac{1}{2} \na_{d(\ln{|\na \Om|^2})} {\overline \Om} , $$
where the first equality uses that $(d (\ln{|\na
\Om|^2}))^{\sharp} \in D$ (see (\ref{lem3-1})), hence $d^J
(\ln{|\na \Om|^2}) = d^I (\ln{|\na \Om|^2})$, and the second
equality follows from (\ref{ddJ+dJd}) and (\ref{naOm-N}).

Since $(d^J (\ln{|\na \Om|^2}))^{\sharp} \in D$ and
$D$ is involutive, we also have
$$\langle dd^J (\ln{|\na \Om|^2}), (\Om - {\overline \Om}) \rangle = 0. $$
Thus, taking the inner product of relation (\ref{db-comparison})
with $1/2(\Om - {\overline \Om})$ and using (\ref{s*-s}) and its
analogue for the almost \ka structure $(g,I)$, we derive
\begin{equation*} \label{ROmbarOm}
\langle R(\Om), {\overline \Om} \rangle = \frac{s}{2} +
\frac{1}{8} |\na {\overline \Om}|^2 ,
\end{equation*}
and the desired form of the ${\overline \Om}$-component of $\rho_0$
follows.

\vspace{0.2cm}

Now we prove the relation (\ref{barlem3-1}). The starting point
for this is the  identity (\ref{bianchi-}) established in Lemma
\ref{Bianchi-lem}. Taking into account that the Ricci tensor is
$J$-invariant and that the anti-self-dual Weyl tensor decomposes
into $W^- = W^-_1 + W^-_2 + W^-_3$ (with respect to the almost
complex structure $I$), relation (\ref{bianchi-}) can be written
as
$$ 0 = \nabla_{JZ} \rho_0 + \frac{1}{6} (d s
\wedge Z^{\flat})^- - 2(\delta W^-_1)_Z - 2(\delta W^-_2)_Z -
2(\delta W^-_3)_Z \; .$$ We consider the ${\overline
\Om}$-component of this relation, taking into account that for
$(\delta W^-_1)_Z$, $(\delta W^-_2)_Z$ we have similar
expressions to the ones for $(\delta W^+_1)_Z$, $(\delta
W^+_2)_Z$ given in (\ref{deltaw^+1}) and (\ref{deltaw^+2}). We thus get
\begin{eqnarray} \nonumber
0 &=& \langle \nabla_{JZ} \rho_0, {\overline \Om} \rangle -
\frac{1}{6}
(d^I ({\bar \kappa} - s))_Z + \frac{1}{2} (d^I {\bar \kappa})_Z \\
\nonumber & & - \langle \na_{IZ} ({\bar \rho}_{*})''_I,
{\overline \Om} \rangle - 2 (\delta ({\bar \rho}_{*})''_I)_Z - 2
\langle (\delta W^-_3)_Z, {\overline \Om} \rangle  \; .
\end{eqnarray}
We plug the expressions for $\rho_0$ and for $({\bar
\rho}_{*})''_I$ given by (\ref{barpsi}) into the above formula to
eventually obtain
\begin{eqnarray} \label{keylem4''}
0 &=& 2(d^I |{\bar a}|^2)_Z -\frac{1}{2} (d^J |{\bar a}|^2)_Z -
\frac{1}{2} (d^I s - d^J s)_Z \\ \nonumber & & - 2 (\delta
{\overline \Psi})_Z - 2 \langle (\delta W^-_3)_Z, {\overline \Om}
\rangle \; .
\end{eqnarray}
To derive the latter formula we have also used that  $ \langle
{\overline
\Psi}, \na_{IZ - JZ} {\overline \Om} \rangle =0$ which follows by Lemma
\ref{step1}, since
$IZ - JZ
\in D^{\perp}$.

Substituting the expression for $\rho_0$ found in (\ref{barpsi})
into (\ref{riccibianchi}), and applying $J$, we get
$$ d^Is- d^Js  = -d^I|{\bar a}|^2 + 4 \delta
{\overline \Psi} . $$ Using this last relation back in (\ref{keylem4''})
to replace the term
$\frac{1}{2} (d^I s - d^J s)_Z$, we obtain
\begin{equation} \label{dbara}
0 = 2(d^I |{\bar a}|^2)_Z -\frac{1}{2} (d^J |{\bar a}|^2 + d^I
|{\bar a}|^2)_Z - 2 \langle (\delta W^-_3)_Z, {\overline \Om}
\rangle \; .
\end{equation}
Since ${\bar a}^{\sharp} \in D$, the 1-form
$$\langle \delta W_3^{-}, {\overline \Om} \rangle =
2W_3^{-}({\bar \phi})({\bar a})$$ is dual to a vector field in
$D^{\perp}$, and then relation (\ref{dbara}) is
equivalent with the relation (\ref{barlem3-1}). $\Box$

\vspace{0.2cm}

The next result is an ``iteration'' of Lemma \ref{step1}.

\begin{Lemma} \label{step3}
 Let $(M,g,J)$ be a
strictly almost \ka 4-manifold satisfying {\rm (\ref{gray3})} and
let $U$ denote the open  set of points where $\na J \neq 0$. Let
$(g,I)$ be the opposite almost \ka structure on $U$, associated to
$(g,J)$ and assume that on a non-empty open subset $U_0 \subset
U$ the opposite almost \ka structure $(g,I)$ is not K\"ahler.
Then, on $U_0$, the opposite almost Hermitian structure $(g,
{\tilde J})$ associated to  $(g,I)$ is $(g, -J)$. Moreover, the
following relation between the canonical Chern forms of $(g,I)$
and $(g,{\tilde J})$ holds
\begin{equation} \label{gammaJ-gammaI}
 \gamma_{\tilde J}= - \gamma_J = 3\gamma_I - dd^I(\ln |\na {\overline
\Om}|^2) -2dd^J(\ln |\na \Om|^2).
\end{equation}
\end{Lemma}

\noindent {\it Proof.} From the assumption that $(g, I)$ is not
K\"ahler on $U_0$ and from Lemma \ref{step1}, it follows that the
\ka nullity of $(g,I)$ is precisely the distribution $D^{\perp}$
on $U_0$; the opposite almost complex structure ${\tilde J}$
associated to  $(g,I)$ is then  defined by ${\tilde J} = I$ on
$D^{\perp}$ and ${\tilde J} = -I$ on $D$, i.e., ${\tilde J} = -
J$.

\vspace{0.2cm}

We now prove (\ref{gammaJ-gammaI}).
As in Lemma \ref{step1}, we need to first investigate
the 2-jet of ${\overline \Om}$ (equivalently, the 1-jet of ${\bar
a}$), with the help of the Ricci relations (\ref{bar*}) and
(\ref{barlem3-1}). Let $p_i, q_i, \; i=1,2$ be the 1-forms defined by the
first relation below. The other three relations follow by (\ref{na-om})
and (\ref{na-phi}).
\begin{equation} \label{nabara} \begin{split}
 \na {\bar a} &=  p_1 \otimes {\bar a} + q_1 \otimes J{\bar a} + p_2 \otimes
J_1 {\bar a} + q_2 \otimes J_2 {\bar a} \; ; \\
 \na (J{\bar a}) &= -q_1
\otimes {\bar a} + p_1 \otimes J{\bar a} + (a-q_2) \otimes J_1
{\bar a} + (p_2-Ja) \otimes J_2 {\bar a} \; ; \\
 \na (J_1 {\bar a}) &=  -p_2
\otimes {\bar a} + (q_2-a) \otimes J{\bar a} + p_1 \otimes J_1
{\bar a} + (b - q_1) \otimes J_2 {\bar a} \; ; \\
 \na (J_2 {\bar a}) &=  -q_2
\otimes {\bar a} + (Ja-p_1) \otimes J{\bar a} + (q_1-b) \otimes
J_1 {\bar a} + p_1 \otimes J_2 {\bar a} \;  . \end{split}
\end{equation}
From the formula (\ref{barlem3-1}) of Lemma \ref{step2}, it
follows immediately that the 1-form $p_1$ is given by
\begin{equation} \label{p1}
p_1 = \frac{1}{2}d(\ln{|\na {\overline \Om}|^2}) = -
\frac{1}{|{\bar a}|^2} I W_3^-({\bar \phi})({\bar a}) .
\end{equation}
An important point of this formula is that $p_1^{\sharp} \in
D^{\perp}$. For the 1-forms $q_1, p_2, q_2$, one obtains
\begin{equation} \label{na-bara} \begin{split}
q_1 &= - Ip_1 - {\bar b} - d^J (\ln{|\na \Om|^2}) ; \\
p_2 &= -\frac{1}{2} I_2 J_1 {\bar a} + \frac{1}{2} Ja ; \\
q_2 &= I p_2 , \end{split}
\end{equation}
where  $I_1$ and
$I_2$ denote the (local) almost complex structures corresponding to
${\bar \phi}$ and $I{\bar \phi}$; recall that $J_1$, $J_2$ stand for
 the almost complex structures dual to $\phi$ and $J\phi$.
Let us just briefly explain how the relations (\ref{na-bara}) are
obtained. Since both $D$ and $D^{\perp}$ are involutive
distributions, we have
\begin{equation*} \begin{split}
(d{\bar a})(X,X') &= (d (J{\bar a}))(X,X') = 0, \;
\forall X, X' \in D^{\perp} \\
 (d (J_1 {\bar a}) )(Y,Y') &= (d (J_2 {\bar a}) )(Y,Y') = 0, \;
\forall \ Y, Y' \in D .
\end{split}
\end{equation*}
and using (\ref{nabara}), these can be seen to be equivalent to
$q_2 = I p_2$. For the remaining two relations we use
(\ref{nabara}) and $q_2 = I p_2$ plugged into the identities
(\ref{bar*}). Taking the $1/2(\Om+ {\overline \Om})$-components
and the $\leftr \La_J^{0,2} M \rightr$-components of the
identities thus obtained, we eventually get
$$ q_1 + I p_1 + {\bar b} = - d^J ( \ln{|\na \Om|^2}), $$
which is the second relation of (\ref{na-bara}); the $\leftr
\La_I^{0,2} M \rightr$-components of the Ricci identities
(\ref{bar*}) imply
$$ p_2^{D} = -\frac{1}{2} I_2 J_1 {\bar a} ,$$
or, equivalently,
$$ p_2 = -\frac{1}{2} I_2 J_1 {\bar a} + p_0 ,
\mbox{ with } p_0^{\sharp} \in D^{\perp} . $$ So far we have
basically followed the track of the computations made in Lemma
\ref{step1}. However, now we can obtain even more by determining
completely the 1-form $p_0$. For this, write the relation between
$\Om$ and ${\overline \Om}$ as
$$ \Om = \frac{2}{|{\bar a}|^2} {\bar a} \wedge J {\bar a} - {\overline
\Om} ,$$ and take the covariant derivative of this relation, using
(\ref{na-om}), (\ref{na-bar-Om}), (\ref{nabara}-\ref{na-bara}).
This eventually leads to $p_0 = \frac{1}{2} Ja $.

\vspace{0.2cm}

As we have already noticed, the opposite almost Hermitian structure $(g,
{\tilde J})$ associated with $(g,I)$ is not a new structure, but  just
$(g, -J)$. However, the iteration of the construction, imposes a natural
new gauge ${\tilde \phi} \in \leftr \La_J^{0,2} M \rightr$,
\begin{equation} \label{tildegauge}
{\tilde \phi} = {\bar \phi} + \frac{2}{|{\bar a}|^2}
I{\bar a} \wedge I{\bar \phi}({\bar a}),
\end{equation}
which may be different than the initial gauge $\phi$. Taking the
covariant derivative of (\ref{tildegauge}) and using
(\ref{na-bar-phi}), (\ref{barlem3-1}) and
(\ref{nabara}-\ref{na-bara}), one finds, after a long
but straightforward computation, that the 1-form ${\tilde b}$ of
the structure $(g,{\tilde J})$ with respect to the gauge ${\tilde
\phi}$ is given by (compare with (\ref{bar-b})):
\begin{equation*} \label{btilde-barb}
{\tilde b} =  3{\bar b} + d^I (\ln |\na {\overline \Om}|^2) + 2
d^J (\ln |\na \Om|^2) .
\end{equation*}
Now relation (\ref{gammaJ-gammaI}) immediately follows since
$d{\tilde b} = - \gamma_{\tilde J} = \gamma_J$. $\Box$

\vspace{0.2cm}
\noindent
{\it Proof of Proposition \ref{localadd-prop}.} We assume as in Lemma
\ref{step3} that $(g,I)$ is {\it not} \ka on an open set $U_0 \subset U$.
Then on
$U_0$ both relations (\ref{gammaI-gammaJ}) and (\ref{gammaJ-gammaI})
hold, and taking the appropriate linear combination of them, we obtain
\begin{equation*} \label{5gamma}
5(\gamma_I - \gamma_J) = dd^I ( \ln{|\na {\overline \Om}|^2} ).
\end{equation*}
On the other hand, starting from (\ref{chernform}) and its
companion relation giving $\gamma_I$, and by using (\ref{s*-s})
and the curvature information picked up in Lemma \ref{step2}, we
compute
\begin{equation*} \label{gammaIJ} \begin{split}
\gamma_I - \gamma_J &= R({\overline \Om}) - {\bar a} \wedge
I{\bar a}
         - R(\Om) + a \wedge Ja \\
  & = - \frac{|\na \Om|^2}{8} (\Om + {\overline \Om})
      - \frac{|\na {\overline \Om}|^2}{16} (\Om - {\overline \Om}).
\end{split}
\end{equation*}
From the last two identities  we obtain
the equality
\begin{equation} \label{keygammaIJ}
 - \frac{|\na \Om|^2}{8} (\Om + {\overline \Om})
      - \frac{|\na {\overline \Om}|^2}{16} (\Om - {\overline \Om}) =
 \frac{1}{5} dd^I ( \ln{|\na {\overline \Om}|^2} ),
\end{equation}
which holds at every point of the open set $U_0$. Since $d^I (
\ln{|\na {\overline \Om}|^2} )$ is dual to a vector field in
$D^{\perp}$ and $D^{\perp}$ is involutive, it follows that the
right hand-side of (\ref{keygammaIJ}) has a zero inner product
with $\Om + {\overline \Om}$. For this to hold for the left
hand-side of (\ref{keygammaIJ}) as well, we must have $|\na \Om| =
0$ at every point of $U_0$. But that is a contradiction, since
$U_0$ is a subset of $U$, the set of points where $|\na \Om| \neq
0$. Therefore, the assumption that the opposite structure $(g, I)$
is not \ka must be false, hence we proved Proposition
\ref{localadd-prop}.  $\Box$

\section{Proof of Theorem \ref{thm2}}

As a consequence of Theorem \ref{thm1} and Proposition \ref{ucp}, we first
obtain the following global version of Proposition \ref{localadd-prop},
which is essential for the proof of Theorem \ref{thm2}, but also
presents interest in its own.
\begin{prop}\label{prop3} Any connected,
strictly almost \ka 4-manifold $(M,g,J,\Om)$ satisfying the third
Gray condition {\rm (\ref{gray3})} admits a globally defined \ka
structure $(g,I)$ which yields the opposite orientation to the
one of $(M,J)$.
\end{prop}
\noindent {\it Proof.} We show that the opposite almost \ka structure
$(g,I)$ defined as in Theorem \ref{thm3} can be extended globally
on $M$. Let $p \in M - U$ be any point of the zero locus of the
Nijenhuis tensor $N$. By the unique continuation property of $N$
(see Proposition \ref{ucp}), there exists a number $k \geq 2$
such that
$$(\na^{\ell} \Om) (p) = 0, \ \forall \ 1\le \ell \le
(k-1); \ \ (\na^k \Om)(p) \neq 0. $$ As a consequence, we obtain
$$0= \langle \na^*\na (\na^{k-1}\Om) , \na^{k-1} \Om \rangle _p =
-\frac{1}{2}\Delta(|\na^{k-1} \Om|^2)(p) + |\na^k\Om|^2(p),$$
which shows that $\Delta (|\na^{k-1} \Om|^2) \neq 0$ at $p$,
hence also in a small neighborhood of $p$. On $U$, the
distribution $D^{\perp}$ is spanned by the commuting hamiltonian
Killing fields  $\frac{\partial}{\partial z}$ and
$\frac{\partial}{\partial t}$ (see Theorem
\ref{thm1}), so that
$d|\na^{k-1}\Om|^2$ is zero on
$D^{\perp}$; equivalently, $|\na^{k-1} \Om|^2$ is a function on the
Riemann surface $\Sigma$. It then follows that  on
$U$ we have (see (\ref{laplacian})):
\begin{equation}\label{extending}
dd^J(|\na^{k-1} \Om|^2) = -\Delta (|\na^{k-1} \Om|^2) \Om_{\Sigma},
\end{equation}
where, recall, $\Om_{\Sigma}$ denotes the volume form of the
(locally defined) Riemann surface $\Sigma$; as a global object
 on
$U$, $\Om_{\Sigma}$ is given by the restriction of
$\Omega$ to $D \subset TU$, i.e. $\Om_{\Sigma}(\cdot,\cdot)= \Om({\rm
pr}^D\cdot, {\rm
pr}^D\cdot)$ where ${\rm pr}^D$ denotes the projection to $D$.  We thus
can define
$I$ in a small neighborhood of
$p$ by setting
$$\Omega_I = -\Om - \frac{2}{\Delta(|\na^{k-1} \Om|^2)} dd^J(|\na^{k-1}
\Om|^2).$$  By (\ref{extending}), this agrees with the definition
of $I$ on $U$. Since $U$ is dense in $M$ and $(g ,I)$ is \ka on
$U$, we conclude that $I$ can be extended as a \ka structure on
whole $M$. $\Box$

\vspace{0.2cm}
Now we are ready to prove Theorem \ref{thm2}.

\vspace{0.2cm} Suppose for contradiction that $(M,g,J,\Om)$ is a
{\it compact} strictly almost \ka 4-manifold whose curvature
satisfies (\ref{gray3}). According to Proposition \ref{prop3}, we
can consider $(M,g,I)$ as a compact \ka surface which also admits
an almost \ka structure $(g,J,\Om)$, compatible with the opposite
orientation to the one of $(M,I)$; equivalently, $\Om$ is an {\it
indefinite} \ka structure on $(M,I)$ (see \cite{petean}).  Note
then that the distributions $D$ and $D^{\perp}$ are well defined
on the whole $M$, respectively as the $(-1)$ and  the
$(+1)$-eigenspace of the symmetric endomorphism $Q= J\circ I$ of
$TM$. Moreover, since $(g,I)$ is K\"ahler, the relation
(\ref{barpsi}) of Lemma \ref{step2} implies that on $U$ (and
therefore, by continuity, on $M$) the Ricci tensor has the form
\begin{equation} \label{ric}
{\rm Ric} = \frac{s}{2}g_\Sigma,
\end{equation}
where $g_\Sigma$ is  defined on whole $M$ as the restriction of $g$
to the distribution $D$. Alternatively, (\ref{ric}) can be directly
derived from the explicit form of $g$ given by (\ref{g0}) of Theorem
\ref{thm1}. Using (\ref{ric}),
it is easy to compute the Chern forms
${\gamma}_I$ and
$\gamma_J$ of
$(g,I)$ and
$(g,J)$, respectively,
\begin{equation}\label{chernforms}
\gamma_I = \frac{s}{2} \Om_{\Sigma}, \ \ \
\gamma_J = (\frac{s}{2} + \frac{|\na \Om|^2}{4})\Om_{\Sigma}.
\end{equation}
 It follows  that $\gamma_I$ and ${\gamma}_J$ are both degenerate on
$M$. From Wu's relations, we conclude
$$c_1^2(M,I) = 2e(M) + 3\sigma(M) = 0; \ c_1^2(M,J) = 2e(M) - 3\sigma(M)
= 0,$$ hence the Euler characteristic $e(M)$ and the signature
$\sigma(M)$ of $(M,I)$ both vanish. From the Kodaira
classification \cite{BPV} and the results of \cite{kotschick0} and
\cite{petean}, we conclude that $(M,I)$ can be a ruled surface
over an elliptic curve, a complex torus, a hyperelliptic surface,
or a minimal properly elliptic surface which is a fibration over a
curve of genus at least two, with no fibers of singular reduction.
We can easily exclude the cases when $(M,I)$ is a torus or a
hyperelliptic surface. Indeed, by replacing $M$ with a finite
cover in the case when $(M,I)$ is hyperelliptic, we can always
assume that $(M,I)$ is a complex torus. By results of Taubes
\cite{taubes1,taubes2}, we know that on a 4-torus any symplectic
form is homotopic to an invariant one, i.e. $[\gamma_I]=0$ and
$[\gamma_J]=0$, so that by (\ref{chernforms}) we obtain
$$0=[\gamma_J - \gamma_I]\cdot [\Omega] = -\frac{1}{4}\int_M  |\na \Om|^2
d\mu,$$ a contradiction with the assumption  that $J$ is not
integrable.

We thus have to analyze the remaining two cases, when $(M,I)$ is a
ruled surface or a minimal elliptic surface. Since
$\frac{\partial}{\partial t}$ (and $\frac{\partial}{\partial z}$)
is a (real) $I$-holomorphic vector field in $D^{\perp}$, we
conclude that $D^{\perp}$ is a holomorphic, complex rank one
distribution on $(U,I)$, hence also on $M$, which gives rise to a
non-singular {\it holomorphic} foliation ${\cal F}$ on the complex
surface $(M,I)$.  From (\ref{chernforms}) we get
\begin{equation}\label{gammaF}
\frac{1}{2}(\gamma_I - \gamma_J) = -\frac{|\na \Om|^2}{8}
\Om_{\Sigma} \in 2\pi c_1(T{\cal F}).
\end{equation}
Our next aim is to show that after replacing $M$ with a finite
cover if necessary, $H^0(M,T{\cal F})\neq 0$.

\vspace{0.2cm}

In the case when $(M,I)$ admits an elliptic fibration (which also
includes some ruled surfaces over an elliptic base), we can use an
argument from \cite{ADK}: since $e(M)=\sigma(M)=0$, we know that
the elliptic fibers are smooth or multiple; since $c_1(M,I) \neq
0$, we have by Kodaira's formula that $c_1(M,I)$ is a non-zero
multiple of the Poincar\'e-dual of any regular (elliptic) fiber
$F$. From (\ref{chernforms}), (\ref{gammaF}) and using the
Poincar\'e duality, we derive
$$0= \frac{1}{2\pi}\int_M \gamma_I\wedge(\gamma_J - \gamma_I)
= c_1(M,I)\cdot [\gamma_J-\gamma_I] = -\frac{1}{8\pi}\int_F |\na
\Om|^2\Omega_{\Sigma},$$ where
$\cdot$ denotes the cup-product of $H^2(M,{\mathbb R})$.  Note
that $\Omega_{\Sigma}$ is a semi-positive (1,1)-form on $(M,I)$
whose kernel, at any point of $U$, is $T{\cal F}$, while $F$ is a
(smooth) holomorphic curve in $(M,I)$; it then follows from the
above equality that, at any point $x\in U$, the fiber $F$ must be
tangent to ${\cal F}$; since $U$ is dense in $M$, we conclude that
${\cal F}$ is  tangent to the fibers everywhere. But as observed
in \cite{ADK}, by replacing $M$ with a finite cover if necessary,
$(M,I)$ then admits a globally defined, non-trivial holomorphic
field tangent to the fibers, which is also a holomorphic section
of $T{\cal F}$.

\vspace{0.2cm} Suppose now that $(M,I)$ is  a ruled surface over
an elliptic base $B$. By the previous argument, we may also assume
that $(M,I)$ does not admit any elliptic fibration. According to
the the classification of non-singular holomorphic foliations
(see \cite[Prop.6]{brunella} and \cite[Sec.3]{mont}), under the
above assumptions for $(M,I)$, the following two cases arise:

\vspace{0.15cm}
\noindent {\it Case 1}: ${\cal F}$ is tangent to the rational fibers;
since $(M,I)$ is not an elliptic fibration, according to \cite{maruyama}
we have $H^0(M,T{\cal F}) \neq 0$ as claimed.

\vspace{0.15cm} \noindent {\it Case 2}:  ${\cal F}$  is a
foliation transversal to the rational fibers, i.e. a {\it Riccati
foliation} in the terminology of \cite{mont}. In this case
$c_1(T{\cal F})$ is the pull-back of a class of $H^2(B,{\mathbb
Z})$ (see \cite{mont}), so that $c_1(T{\cal F})$ is a multiple of
the Poincar\'e-dual of any rational fiber $F$. Since $[F]\cdot
[F] = 0$ in homology and in view of (\ref{gammaF}), we then have
$$0 = \frac{1}{2} \int_F (\gamma_I - \gamma_J) = -\frac{1}{16\pi}\int_F |\na
\Om|^2\Omega_{\Sigma}. $$ Then, we conclude as in the case of
elliptic fibrations that each fiber $F$ must be tangent to ${\cal
F}$, which is a contradiction.

\vspace{0.2cm}

As a final step of the proof of Theorem \ref{thm2}, note that
$T{\cal F}= D^{\perp}$ lies in the kernel of the Ricci tensor of
$g$ (see (\ref{ric})); then, by the well-known
Bochner-Lichnerowicz argument, for any holomorphic section $\Xi$
of $T{\cal F} \subset TM$ we get $\na \Xi =0$. Since we already
showed $H^0(M,T{\cal F}) \neq 0$, it follows that $T{\cal F}$
should be parallel, and then, so must be its orthogonal complement
$D$. But then $J$ must be parallel too, a contradiction. $\Box$

\vspace{0.2cm} \noindent {\bf Acknowledgments.} The authors are
grateful to O. Biquard, D. Calderbank and P. Gauduchon for their
interest and stimulating discussions.


\begin{thebibliography}{99}

\bibitem{AGI} B. Alexandrov, G. Grantcharov and S. Ivanov, {\it Curvature
properties of twistor spaces of quaternionic K\"ahler manifolds}, J. Geom.
{\bf 62} (1998), 1--12.

\bibitem{AA} V. Apostolov and J. Armstrong, {\it Symplectic 4-manifolds
with Hermitian Weyl tensor}, Trans. Amer. Math. Soc. {\bf 352} (2000),
4501--4513.

\bibitem{AAD} V. Apostolov, J. Armstrong and  T. Draghici, {\it Local
rigidity of certain classes of almost K\"ahler 4-manifolds}, to
appear in Ann. Glob. Anal. Geom.

\bibitem{ACG} V. Apostolov, D. Calderbank and P. Gauduchon, {\it The
geometry of weakly selfdual K\"ahler surfaces}, to
appear in Compositio Math.


\bibitem{ADK} V. Apostolov, T. Draghici and D. Kotschick, {\it An
integrability theorem for almost K{\"a}hler 4-manifolds}, C. R. Acad.
Sci. Paris {\bf 329}, s{\'e}r. {I} (1999), 413--418.

\bibitem{ADM} V. Apostolov, T. Draghici and A. Moroianu, {\it A
splitting theorem for K\"ahler manifolds whose Ricci tensors have
constant eigenvalues}, Int. J. Math. {\bf 12} (2001), 769--789.


\bibitem{arons} N. Aronszajn, {\it A unique continuation theorem
for solutions of elliptic partial differential equations or
inequalities of second order}, J. Math. Pures Appl. {\bf 36} (1957),
235--249.

\bibitem{Arm0} J. Armstrong, {\it On four-dimensional almost K{\"a}hler
manifolds}, Quart. J. Math. Oxford Ser.(2) {\bf 48} (1997),
405--415.

\bibitem{Arm1} J. Armstrong, {\it An ansatz for Almost-K{\"a}hler,
Einstein 4-manifolds}, J. reine angew. Math. {\bf 542} (2002), 53--84.


\bibitem{BPV} W.Barth, C.Peters and A.Van de Ven,  Compact complex
surfaces, Springer-Verlagh, Berlin  Heidelberg New York Tokyo, 1984.

\bibitem{besse} A. L. Besse, Einstein manifolds, Ergeb. Math. Grenzgeb.,
Springer-Verlag, Berlin, Heidelberg, New York, 1987.


\bibitem{Bl} D. E. Blair, {\it The "total scalar curvature" as a
symplectic
invariant and
related results}, Proc. 3rd Congress of Geometry, Thessaloniki (1991),
79--83.

\bibitem{BI} D. E. Blair and S. Ianus, {\it Critical associated metrics on
symplectic manifolds},
Contemp. Math. {\bf 51} (1986), 23--29.


\bibitem{bourg} J.-P. Bourguignon, {\it Les vari{\'e}t{\'e}s de dimension 4
{\`a} signature non nule dont la courbure est harmonique sont
d'Einstein}, Invent. Math. {\bf 63} (1981), 263--286.

\bibitem{brunella} M. Brunella, {\it Feuilletages holomorphes sur les
surfaces complexes compactes}, Ann. Scient. \'Ec. Norm. Sup., S\'er. 4
{\bf 30} (1997), 569--594.

\bibitem{DM} J. Davidov and O. Mu\u{s}karov, {\it Twistor spaces with
Hermitian
Ricci
tensor}, Proc. Amer. Math. Soc. {\bf 109} (1990), 1115--1120.

\bibitem{donaldson1} S. K. Donaldson, {\it Remarks on Gauge theory,
complex geometry and 4-manifolds topology}, in ``The Fields Medallists
Lectures'' (eds. M. Atiyah and D. Iagolnitzer), pp. 384--403, World
Scientific, 1997.

\bibitem{donaldson2} S. K. Donaldson, {\it Symmetric spaces, K{\"a}hler
geometry and Hamiltonian dynamics}, in ``Northern California Symplectic
Geometry Seminar'' (eds. Y. Eliashberg et al.), pp.13--33, American
Mathematical Society, 1999.

\bibitem{draghici0} T. Draghici,
{\it On some $4$-dimensional almost K{\"a}hler manifolds}, Kodai
Math. J. {\bf 18} (1995), 156--163.


\bibitem{Dr2} T. Draghici, {\it Almost \ka 4-manifolds with
$J$-invariant Ricci tensor}, Houston J. Math. {\bf 25} (1999),
133--145.


\bibitem{gauduchon1} P. Gauduchon, {\it Hermitian connections and Dirac
operators}, Boll. U.M.I. (7) {\bf 11}-B (1997), Suppl. fasc. 2,
257--288.

\bibitem{GH} G. Gibbons and  S. Hawking, {\it Classification of
Gravitational Instanton Symmetries}, Comm. Math. Phys. {\bf 66} (1979),
291--310.

\bibitem{Go} S. I. Goldberg, {\it Integrability of almost \ka manifolds},
Proc. Amer. Math. Soc. {\bf 21} (1969), 96--100.

\bibitem{mont} X. G\'omez-Mont, {\it Holomorphic foliations in ruled
surfaces}, Trans. Amer. Math. Soc. {\bf 312} (1989), 179--201.

\bibitem{gray} A. Gray, {\it Riemannian manifolds with geodesic symmetries of
order $3$}. J. Differential Geom. {\bf 7}
(1972), 343--369.

\bibitem{Gr} A. Gray, {\it Curvature identities for Hermitian and almost
Hermitian manifolds}, T{\^o}hoku Math. J. {\bf 28} (1976), 601--612.

\bibitem{gromov} M. Gromov, {\it Pseudoholomorphic curves in symplectic
manifolds}, Invent. Math. {\bf 82} (1985), 307--347.

\bibitem{joyce} D. Joyce, {\it Explicit construction of self-dual
4-manifolds}, Duke Math. J. {\bf 77} (1995), 519--552.

\bibitem{kazdan} J. Kazdan, {\it Unique Continuation in Geometry}, Comm.
Pure and Appl. Math. {\bf XLI} (1988), 667--681.

\bibitem{KN2} S. Kobayashi and K. Nomizu, Foundations of differential
geometry, Vol II Interscience Publishers (John Wiley Inc.) No {\bf 15},
1969.

\bibitem{kotschick0} D. Kotschick, {\it Orientations and geometrisations
of compact complex surfaces}, Bull. London Math. Soc. {\bf 29} (1997),
145--149.

\bibitem{kotschick} D. Kotschick, {\it Einstein metrics and smooth
structures}, Geom. Topol. {\bf 2} (1998), 1--10.

\bibitem{kow} O. Kowalski, Generalized Symmetric Spaces, LNM {\bf 805},
1980.

\bibitem{Le-Wang} H.-V. L\^e and  G. Wang, {\it Anti-complexified Ricci
flow on compact symplectic manifolds}, J. reine angew. Math. {\bf 530}
(2001), 17--31.

\bibitem{LeB} C.LeBrun, {\it Explicit self-dual metrics on ${\Bbb
CP}^2 \sharp ... \sharp {\Bbb CP}^2$}, J.Differential
Geom. {\bf 34} (1991), 223--253.


\bibitem{lebrun1} C.LeBrun, {\it Four-manifolds without Einstein
metrics}, Math. Res. Lett. {\bf 3} (1996), 133--147.

\bibitem{lebrun2} C. LeBrun, {\it Weyl curvature, Einstein metrics, and
Seiberg-Witten theory}, Math. Res. Lett. {\bf 5} (1998), 423--438.

\bibitem{lebrun3} C. LeBrun, {\it Ricci curvature, minimal volumes, and
Seiberg-Witten theory}, Invent. Math. {\bf 145} (2001), 279--316.

\bibitem{maruyama} M. Maruyama, {\it On automorphism groups of ruled
surfaces}, J. Math. Koyoto Univ. {\bf 11} (1971), 89--112.

\bibitem{Ne-Ni} A. Newlander and L. Nirenberg, {\it Complex analytic
coordinates in almost complex manifolds}, Ann. Math. {\bf 65} (1957),
391--404.


\bibitem{NuP} P. Nurowski and  M. Przanowski, {\it A four-dimensional
example of Ricci flat metric admitting almost K{\"a}hler non-K{\"a}hler
structure}, Classical Quantum Gravity {\bf 16} (1999), no {\bf 3},
L9--L13.

\bibitem{petean} J. Petean, {\it Indefinite K{\"a}hler-Einstein metrics on
compact complex surfaces}, Comm. Math. Phys. {\bf 189} (1997), 227--235.

\bibitem{Tod} K. P. Tod, private communication to the second author.

\bibitem{Se2} K. Sekigawa, {\it On some compact Einstein almost \ka
manifolds},
J. Math. Soc. Japan {\bf 36} (1987), 677--684.

\bibitem{taubes1} C.H. Taubes, {\it The Seiberg-Witten Invariants and
Symplectic Forms}, Math. Res. Lett. {\bf 1} (1994), 809--822.

\bibitem{taubes2} C.H. Taubes, {\it The Seiberg-Witten and Gromov
Invariants}, Math. Res. Lett. {\bf 2} (1995), 221--238.

\bibitem{TV} F. Tricerri and L. Vanhecke, {\it Curvature tensors on almost
Hermitian manifolds},
Trans. Amer. Math. Soc. {\bf 267} (1981), 365--398.

\bibitem{wood} C. Wood, {\it Harmonic almost-complex structures},
Compositio Math. {\bf 99} (1995), 183--212.

\end{thebibliography}
\end{document}